\documentclass[reqno,10pt,a4paper]{amsart}
\usepackage[margin=2.7cm]{geometry}
\usepackage{cite}
\usepackage{newpxtext}
\usepackage{cabin}
\usepackage[varqu,varl]{inconsolata}
\usepackage[bigdelims,vvarbb]{newpxmath}
\usepackage[cal=cm,bb=esstix,bbscaled=1.05]{mathalfa}
\usepackage{graphicx}
\usepackage{colortbl}
\usepackage{algorithm}
\usepackage{algpseudocode}
\usepackage{appendix}
\usepackage{booktabs}
\usepackage[up]{caption}
\usepackage[labelformat=simple]{subcaption}
\usepackage{xcolor}
\usepackage{hyperref}
\hypersetup{
    colorlinks,
    linkcolor={red!50!black},
    citecolor={blue!50!black},
    urlcolor={blue!80!black}
}

\numberwithin{equation}{section}
\AtBeginDocument{%
\def\MR#1{}
}
\newcommand\blfootnote[1]{%
  \begingroup
  \renewcommand\thefootnote{}\footnote{#1}%
  \addtocounter{footnote}{-1}%
  \endgroup
}
\newcommand{\orcid}[1]{\,\resizebox{8px}{!}{\href{https://orcid.org/#1}{\includegraphics{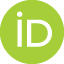}}}}
\graphicspath{{images_300/}} 
\newcommand{\ZZ}{\mathbb{Z}}
\newcommand{\PP}{\mathbb{P}}
\newcommand{\CC}{\mathbb{C}}
\newcommand{\QQ}{\mathbb{Q}}
\newcommand{\RR}{\mathbb{R}}

\newcommand{\Cstar}{\CC^\times}

\newcommand{\QFano}{$\QQ$\nobreakdash-Fano}
\newcommand{\Qfactorial}{$\QQ$\nobreakdash-factorial}
\DeclareMathOperator{\SL}{SL}
\DeclareMathOperator{\Cone}{cone}
\DeclareMathOperator{\conv}{conv}
\newcommand{\doi}[1]{\href{https://doi.org/#1}{doi:#1}}

\theoremstyle{plain} 
\newtheorem{Lemma}{Lemma}
\newtheorem{Thm}[Lemma]{Theorem}

\newtheorem{Prop}[Lemma]{Proposition}

\theoremstyle{definition}
\newtheorem{Defn}[Lemma]{Definition}

\newtheorem{Not}[Lemma]{Notation}
\theoremstyle{remark}

\begin{document}
\author[T.\,Coates]{Tom Coates\orcid{0000-0003-0779-9735}}
\address{Department of Mathematics\\Imperial College London\\180 Queen's Gate\\London\\SW7 2AZ\\UK}
\email{t.coates@imperial.ac.uk}
\author[A.\,M.\,Kasprzyk]{Alexander M.\ Kasprzyk\orcid{0000-0003-2340-5257}}
\address{School of Mathematical Sciences\\University of Nottingham\\Nottingham\\NG7 2RD\\UK}
\email{a.m.kasprzyk@nottingham.ac.uk}
\author[S.\,Veneziale]{Sara Veneziale\orcid{0000-0003-2851-3820}}
\address{Department of Mathematics\\Imperial College London\\180 Queen's Gate\\London\\SW7 2AZ\\UK}
\email{s.veneziale21@imperial.ac.uk}
\makeatletter
\@namedef{subjclassname@2020}{%
  \textup{2020} Mathematics Subject Classification}
\makeatother
\title{Machine learning detects terminal singularities}
\begin{abstract}
Algebraic varieties are the geometric shapes defined by systems of polynomial equations; they are ubiquitous across mathematics and science. Amongst these algebraic varieties are \QFano{} varieties: positively curved shapes which have \Qfactorial{} terminal singularities. \QFano{} varieties are of fundamental importance in geometry as they are `atomic pieces' of more complex shapes -- the process of breaking a shape into simpler pieces in this sense is called the Minimal Model Programme.

Despite their importance, the classification of \QFano{} varieties remains unknown. In this paper we demonstrate that machine learning can be used to understand this classification. We focus on eight-dimensional positively-curved algebraic varieties that have toric symmetry and Picard rank two, and develop a neural network classifier that predicts with~95\% accuracy whether or not such an algebraic variety is \QFano. We use this to give a first sketch of the landscape of \QFano{} varieties in dimension eight.

How the neural network is able to detect \QFano{} varieties with such accuracy remains mysterious, and hints at some deep mathematical theory waiting to be uncovered. Furthermore, when visualised using the quantum period, an invariant that has played an important role in recent theoretical developments, we observe that the classification as revealed by~ML appears to fall within a bounded region, and is stratified by the Fano index. This suggests that it may be possible to state and prove conjectures on completeness in the future.

Inspired by the~ML analysis, we formulate and prove a new global combinatorial criterion for a positively curved toric variety of Picard rank two to have terminal singularities. Together with the first sketch of the landscape of \QFano{} varieties in higher dimensions, this gives strong new evidence that machine learning can be an essential tool in developing mathematical conjectures and accelerating theoretical discovery.
\end{abstract}
\maketitle
\blfootnote{37th Conference on Neural Information Processing Systems (NeurIPS 2023).}
\section{Introduction}\label{sec:introduction}
\begin{figure}[bhtp]
    \centering
    \begin{subfigure}{.2\textwidth}
      \includegraphics[width=\linewidth]{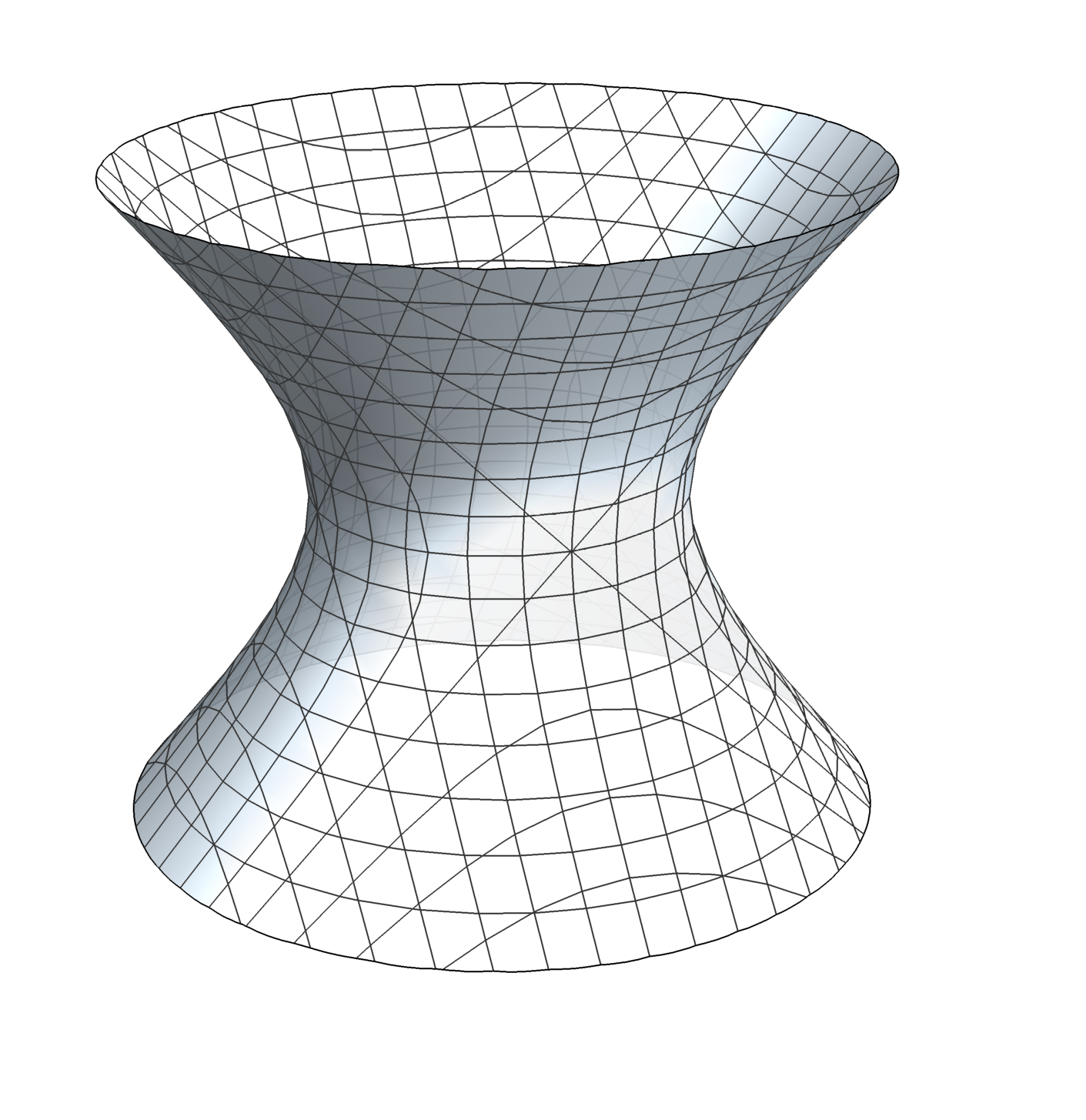}
      \caption{\ $x^2+y^2=z^2+1$}
      \label{fig:non_singular}
    \end{subfigure}
    \qquad \quad 
    \begin{subfigure}{.2\textwidth}
      \includegraphics[width=\linewidth]{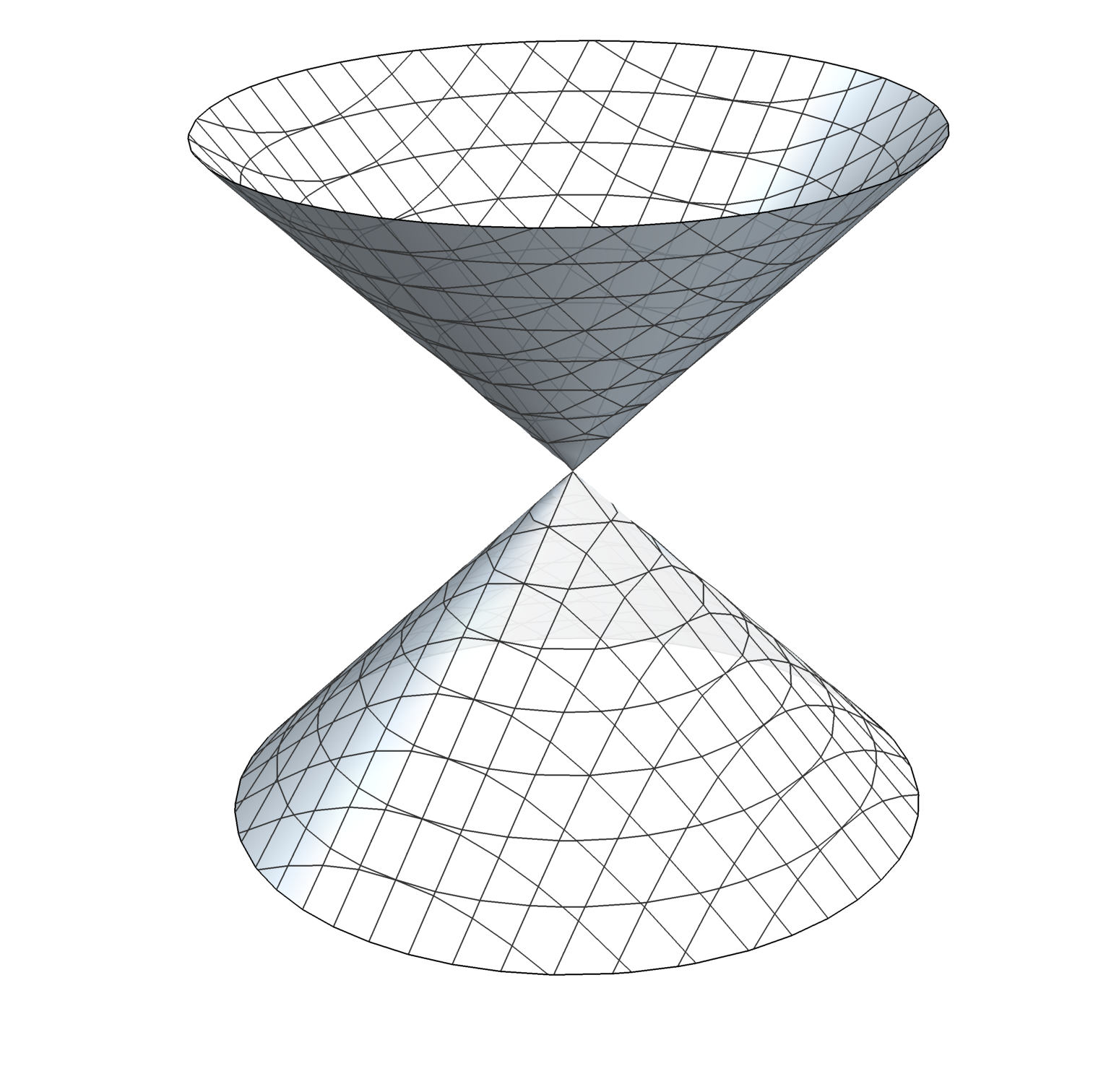}
      \caption{\ $x^2+y^2 = z^2$}
      \label{fig:singular}
    \end{subfigure}
    \qquad \quad 
    \begin{subfigure}{.2\textwidth}
        \includegraphics[width=\linewidth]{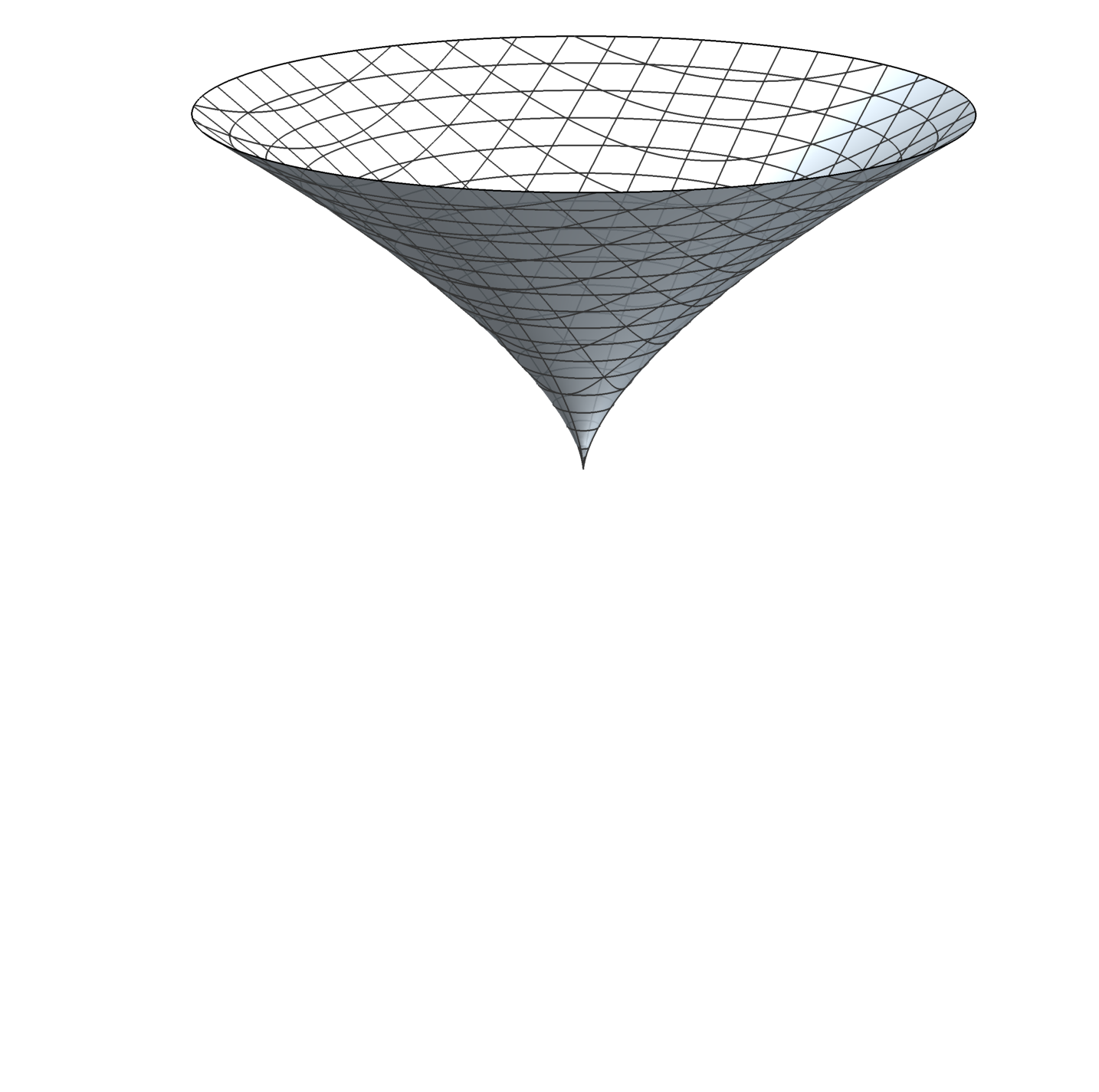}
        \caption{\ $x^2+y^2 = z^3$}
        \label{fig:singular2}
      \end{subfigure}
    \caption{Algebraic varieties in $\RR^3$ with different defining equations.}
    \label{fig:singular_or_not}
  \end{figure}

Systems of polynomial equations occur throughout mathematics and science; see e.g.~\cite{AtiyahDrinfeldHitchinManin1978,Greene1997,ErikssonRanestadSturmfelsSullivant2005,NiederreiterXing2009}. Solutions of these systems define shapes called~\emph{algebraic varieties}. Depending on the equations involved, algebraic varieties can be smooth (as in Figure~\ref{fig:non_singular}) or have singularities (as in Figures~\ref{fig:singular} and~\ref{fig:singular2}). In this paper we show that machine learning methods can detect a class of singularities called~\emph{terminal singularities}.

A key class of algebraic varieties are~\emph{Fano varieties}:  positively curved shapes that are basic building blocks in algebraic geometry. Fano varieties are `atomic pieces' of more complex shapes, in the sense of the Minimal Model Programme~\cite{Kollar1987,KollarMori1998,Cascini2021}. Running the Minimal Model Programme -- that is, breaking an algebraic variety~$X$ into atomic pieces -- involves making birational transformations of~$X$. These are modifications on subsets with zero volume (and codimension at least one), and can either introduce or remove singularities. The building blocks that emerge from this process are not necessarily smooth: they satisfy a weaker condition called \Qfactorial{}ity,\footnote{An algebraic variety~$X$ is \Qfactorial{} if it is normal and, in addition, for each rank-one reflexive sheaf~$E$ on~$X$, some tensor power of~$E$ is a line bundle. This implies that the dimension of the singular locus in~$X$ is at most~$\dim X-2$, and that some tensor power of the canonical sheaf (of top-degree differential forms) is a line bundle.} and can have mild singularities called terminal singularities~\cite{Reid1987}. Fano varieties that are \Qfactorial{} and have terminal singularities are called~\emph{\QFano{} varieties}.

The classification of \QFano{} varieties is therefore a long-standing problem of great importance~\cite{DelPezzo1887,MoriMukai1981,MoriMukai1981:erratum,KollarMiyaokaMoriTakagi2000,Birkar2021} -- one can think of this as building a Periodic Table for geometry. But, despite more than a century of study, very little is known. In what follows we exploit the fact that machine learning can detect terminal singularities to give the first sketch of part of the classification of higher\nobreakdash-dimensional \QFano{} varieties.

We probe the classification of \QFano{} varieties using a class of highly-symmetrical shapes called~\emph{toric varieties}. (For example, the algebraic varieties pictured in Figure~\ref{fig:singular_or_not} are toric varieties.) Toric varieties are particularly suitable for computation and machine learning, because their geometric properties are encoded by simple combinatorial objects. We consider Fano toric varieties of Picard rank two. These can be encoded using a~$2 \times N$ matrix of non-negative integers called the~\emph{weight matrix}; here the dimension of the toric variety is~$N-2$. 

To determine whether such a toric variety~$X$ is a \QFano{} variety we need to check whether~$X$ is \Qfactorial{}, and whether the singularities of~$X$ are terminal. Checking \Qfactorial{}ity from the weight matrix of~$X$ turns out to be straightforward (see~\S\ref{sec:data_generation}) but checking terminality is extremely challenging. This is because there is no satisfactory theoretical understanding of the problem. We lack a global criterion for detecting terminality in terms of weight data (such as~\cite{Kasprzyk2013} in a simpler setting) and so have to fall back on first enumerating all the singularities to analyse, and then checking terminality for each singularity. Each step is a challenging problem in discrete geometry: the first step involves building a different combinatorial object associated to the $n$\nobreakdash-dimensional toric variety~$X$, which is a collection of cones in~$\RR^n$ called the~\emph{fan} $\Sigma(X)$; the second step involves checking for various cones in the fan whether or not they contain lattice points on or below a certain hyperplane. To give a sense of the difficulty of the computations involved, generating and post-processing our dataset of 10~million toric varieties in dimension eight took around 30~CPU years.

To overcome this difficulty, and hence to begin to investigate the classification of \QFano{} varieties in dimension eight, we used supervised machine learning. We trained a feed-forward neural network classifier on a balanced dataset of 5~million examples; these are eight\nobreakdash-dimensional \Qfactorial{} Fano toric varieties of Picard rank two, of which 2.5~million are terminal and 2.5~million non-terminal. Testing on a further balanced dataset of 5 million examples showed that the neural network classifies such toric varieties as terminal or non-terminal with an accuracy of~95\%. This high accuracy allowed us to rapidly generate many additional examples that are with high probability \QFano{} varieties -- that is, examples that the classifier predicts have terminal singularities. This ML\nobreakdash-assisted generation step is much more efficient: generating 100~million examples in dimension eight took less than 120~CPU hours.

The fact that the~ML classifier can detect terminal singularities with such high accuracy suggests that there is new mathematics waiting to be discovered here -- there should be a simple criterion in terms of the weight matrix to determine whether or not a toric variety~$X$ has terminal singularities. In~\S\ref{sec:algorithm} we take the first steps in this direction, giving in Algorithm~\ref{alg:terminality_algorithm} a new method to check terminality directly from the weight matrix, for toric varieties of Picard rank two. A proof of correctness is given in~\S\ref{sec:proof}. This new algorithm is fifteen times faster than the na\"ive approach that we used to generate our labelled dataset, but still several orders of magnitude slower than the neural network classifier. We believe that this is not the end of the story, and that the~ML results suggest that a simpler criterion exists. Note that the neural network classifier cannot be doing anything analogous to Algorithm~\ref{alg:terminality_algorithm}: the algorithm relies on divisibility relations between entries of the weight matrix~(GCDs etc.) that are not visible to the neural network, as they are destroyed by the rescaling and standardisation that is applied to the weights before they are fed to the classifier.

In~\S\ref{sec:application} we use the~ML-assisted dataset of 100~million examples to begin to explore the classification of \QFano{} varieties in dimension eight. We visualise the dataset using the~\emph{regularized quantum period}, an invariant that has played an important role in recent theoretical work on \QFano{} classification, discovering that an appropriate projection of the data appears to fill out a wedge-shaped region bounded by two straight lines. This visualisation suggests some simple patterns in the classification: for example, the distance from one edge of the wedge appears to be determined by the Fano index of the variety.

Our work is further evidence that machine learning can be an indispensable tool for generating and guiding mathematical understanding. The neural network classifier led directly to Algorithm~\ref{alg:terminality_algorithm}, a new theoretical result, by revealing that the classification problem was tractable and thus there was probably new mathematics waiting to be found. This is part of a new wave of application of artificial intelligence to pure mathematics~\cite{DaviesEtAl2021, ErbinFinotello2021, Wagner2021, He2022, WuDeLoera2022, CoatesKasprzykVeneziale2022, Williamson2023}, where machine learning methods drive theorem discovery. 

A genuinely novel contribution here, though, is the use of machine learning for data generation and data exploration in pure mathematics. Sketching the landscape of higher\nobreakdash-dimensional \QFano{} varieties using traditional methods would be impossible with the current theoretical understanding, and prohibitively expensive using the current exact algorithms. Training a neural network classifier however, allows us to explore this landscape easily -- a landscape that is unreachable with current mathematical tools.

\subsection*{Why dimension eight?}
We chose to work with eight\nobreakdash-dimensional varieties for several reasons. It is important to distance ourselves from the surface case (dimension two), where terminality is a trivial condition. A two\nobreakdash-dimensional algebraic variety has terminal singularities if and only if it is smooth. On the other hand, we should consider a dimension where we can generate a sufficient amount of data for machine learning (the analogue of our dataset in dimension three, for example, contains only 34~examples~\cite{Kasprzyk2006}) and where we can generate enough data to meaningfully probe the classification. Moreover, we work in Picard rank two because there already exists a fast combinatorial formula to check terminality in rank one~\cite{Kasprzyk2013}; Picard rank two is the next natural case to consider.

\section{Mathematical background}\label{sec:mathematical_background_main}
The prototypical example of a Fano variety is projective space~$\PP^{N-1}$, which can be thought of as the quotient of~$\CC^N \setminus \{\mathbf{0}\}$ by~$\CC^{\times}$ acting as follows:
\[
    \lambda \cdot (z_1, \dots, z_N) = (\lambda z_1, \dots, \lambda z_N) 
\]
Fano toric varieties of Picard rank two arise similarly. They can be constructed as the quotient of~$\CC^N \setminus S$, where~$S$ is a union of subspaces, by an action of~$(\CC^{\times})^2$. This action, and the union of subspaces~$S$, is encoded by a weight matrix:
\begin{align}\label{eq:weight_matrix}
    \begin{bmatrix}
        a_1 & \cdots & a_{N} \\
        b_1 & \cdots & b_{N}
    \end{bmatrix} 
\end{align}
Here we assume that all~$(a_i,b_i) \in \ZZ^2 \setminus \{\mathbf{0}\}$ lie in a strictly convex cone~$C \subset \RR^2$. The action is
\[
    (\lambda, \mu) \cdot (z_1, \dots, z_{N}) = (\lambda^{a_1} \mu^{b_1} z_1 , \dots, \lambda^{a_{N}} \mu^{b_{N}} z_{N} )
\] 
and~$S = S_+ \cup S_-$ is the union of subspaces~$S_+$ and~$S_-$, where
\begin{equation}
    \label{eq:S_plus_minus}
    \begin{aligned}
        &S_+ = \{(z_1, \dots, z_{N}) \mid z_i =0 \text{ if } b_i/a_i > b/a\} \\
        &S_- = \{(z_1, \dots, z_{N}) \mid z_i =0 \text{ if } b_i/a_i < b/a\} 
    \end{aligned}
\end{equation}
and~$a = \sum_{i=1}^{N} a_i$,~$b = \sum_{i=1}^{N} b_i$: see~\cite{BrownCortiZucconi2004}. The quotient~$X= (\CC^N \setminus S)/(\CC^{\times})^2$ is an algebraic variety of dimension~$N-2$. We assume in addition that both~$S_+$ and~$S_-$ have dimension at least two; this implies that the second Betti number of~$X$ is two, that is,~$X$ has Picard rank two.

Since we have insisted that all columns~$(a_i,b_i)$ lie in a strictly convex cone~$C$, we can always permute columns and apply an~$\SL_2(\ZZ)$ transformation to the weight matrix to obtain a matrix in standard form:
\begin{align}\label{eq:standard_form}
    \begin{bmatrix}
        a_1 & a_2 & \cdots & a_{N} \\
        0 & b_2 & \cdots & b_{N}    
   \end{bmatrix}
\end{align}
where all entries are non-negative, the columns are cyclically ordered anticlockwise, and~$a_{N}<b_{N}$. This transformation corresponds to renumbering the co-ordinates of~$\CC^N$ and reparametrising the torus~$(\Cstar)^2$ that acts, and consequently leaves the quotient variety~$X$ that we construct unchanged.

We will consider weight matrices~\eqref{eq:weight_matrix} that satisfy an additional condition called being~\emph{well-formed}. An~$r \times N$ weight matrix is called standard if the greatest common divisor of its~$r \times r$ minors is one, and is well-formed if every submatrix formed by deleting a column is standard~\cite{Ahmadinezhad2017}. Considering only well-formed weight matrices guarantees that a toric variety determines and is determined by its weight matrix, uniquely up to $\SL_r(\ZZ)$\nobreakdash-transformation.

\subsection*{Testing terminality}
As mentioned in the introduction, an $n$\nobreakdash-dimensional toric variety~$X$ determines a collection~$\Sigma(X)$ of cones in~$\RR^n$ called the fan of~$X$. A toric variety is completely determined by its fan. The process of determining the fan~$\Sigma(X)$ from the weight matrix~\eqref{eq:weight_matrix} is explained in~\S\ref{sec:mathematical_background_supp}; this is a challenging combinatorial calculation. In the fan~$\Sigma(X)$, the one\nobreakdash-dimensional cones are called rays. For a Fano toric variety~$X$, taking the convex hull of the first lattice point on each ray defines a convex polytope~$P$, and~$X$ has terminal singularities if and only if the only lattice points in~$P$ are the origin and the vertices. Verifying this is a conceptually straightforward but computationally challenging calculation in integer linear programming.

\section{Data generation}\label{sec:data_generation}
We generated a balanced, labelled dataset of ten million \Qfactorial{} Fano toric varieties of Picard rank two and dimension eight. These varieties are encoded, as described above, by weight matrices. We generated~$2 \times 10$ integer-valued matrices in standard form, as in~\eqref{eq:standard_form}, with entries chosen uniformly at random from the set~$\{0,\ldots,7\}$. Minor exceptions to this were the values for~$a_1$ and~$b_{N}$, which were both chosen uniformly at random from the set~$\{1,\ldots,7\}$, and the value for~$a_{N}$, which was chosen uniformly at random from the set~$\{0,\ldots,b_{N}-1\}$. Once a random weight matrix was generated, we retained it only if it satisfied:
\begin{enumerate}
    \item~\label{item:zeros} None of the columns are the zero vector.
    \item~\label{item:simplicial} The sum of the columns is not a multiple of any of them.
    \item~\label{item:rays} The subspaces~$S_+$ and~$S_-$ in~\eqref{eq:S_plus_minus} are both of dimension at least two.
    \item~\label{item:well-formed} The matrix is well-formed.
\end{enumerate}
The first condition here was part of our definition of weight matrix; the second condition is equivalent to~$X$ being \Qfactorial{}; the third condition guarantees that~$X$ has Picard rank two; and the fourth condition was discussed above.

We used rejection sampling to ensure that the dataset contains an equal number of terminal and non-terminal examples. Before generating any weight matrix, a boolean value was set to \texttt{True} (terminal) or \texttt{False} (non-terminal). Once a random weight matrix that satisfied conditions~\eqref{item:zeros}--\eqref{item:well-formed} above was generated, we checked if the corresponding toric variety was terminal using the method discussed in~\S\ref{sec:mathematical_background_main}. If the terminality check agreed with the chosen boolean, the weight matrix was added to our dataset; otherwise the generation step was repeated until a match was found.

As discussed, different weight matrices can give rise to the same toric variety. Up to isomorphism, however, a toric variety~$X$ is determined by the isomorphism class of its fan. We deduplicated our dataset by placing the corresponding fan~$\Sigma(X)$, which we had already computed in order to test for terminality, in normal form~\cite{GrinisKasprzyk2013,KreuzerSkarke2004}. In practice, very few duplicates occurred.

\section{Building the machine learning model}\label{sec:apply_ml}
We built a neural network classifier to determine whether a \Qfactorial{} Fano variety of Picard rank two and dimension eight is terminal. The network was trained on the features given by concatenating the two rows of a weight matrix,~$[a_1, \dots, a_{10}, b_1, \dots, b_{10}]$. The features were standardised by translating their mean to zero and scaling to variance one. The network, a multilayer perceptron, is a fully connected feedforward neural network with three hidden layers and leaky ReLu activation function. It was trained on the dataset described in~\S\ref{sec:data_generation} using binary cross-entropy as loss function, stochastic mini-batch gradient descent optimiser and using early-stopping, for a maximum of 150~epochs and with learning rate reduction on plateaux. We tested the model on a balanced subset of~50\% of the data~(5M); the remainder was used for training~(40\%; 4M balanced) and validation~(10\%; 1M).

Hyperparameter tuning was partly carried out using RayTune~\cite{Raytune} on a small portion of the training data, via random grid search with Async Successive Halving Algorithm (ASHA) scheduler~\cite{Li2020}, for 100~experiments. Given the best configuration resulting from the random grid search, we then manually explored nearby configurations and took the best performing one. The final best network configuration is summarised in Table~\ref{tab:network_architecture}.

\begin{table}[t]
    \centering
    \vspace{1em}
    \begin{tabular}{rlp{0.01ex}rl}
        \toprule
        \multicolumn{1}{c}{\textbf{Hyperparameter}} &
        \multicolumn{1}{c}{\textbf{Value}} && 
        \multicolumn{1}{c}{\textbf{Hyperparameter}} & 
        \multicolumn{1}{c}{\textbf{Value}} \\
        \cmidrule{1-2} \cmidrule{4-5}
        \texttt{Layers} & $(512,768,512)$  &&
        \texttt{Momentum} & $0.99$\\
        \texttt{Batch size} & $128$ &&
        \texttt{LeakyRelu slope} & $0.01$ \\
        \texttt{Initial learning rate} & $0.01$ \\
        \bottomrule
    \end{tabular}
    \caption{Final network architecture and configuration.}\label{tab:network_architecture}
\end{table}    

\begin{figure}
    \centering
    \begin{subfigure}[t]{0.36\textwidth}
        \centering
        \includegraphics[width=\textwidth]{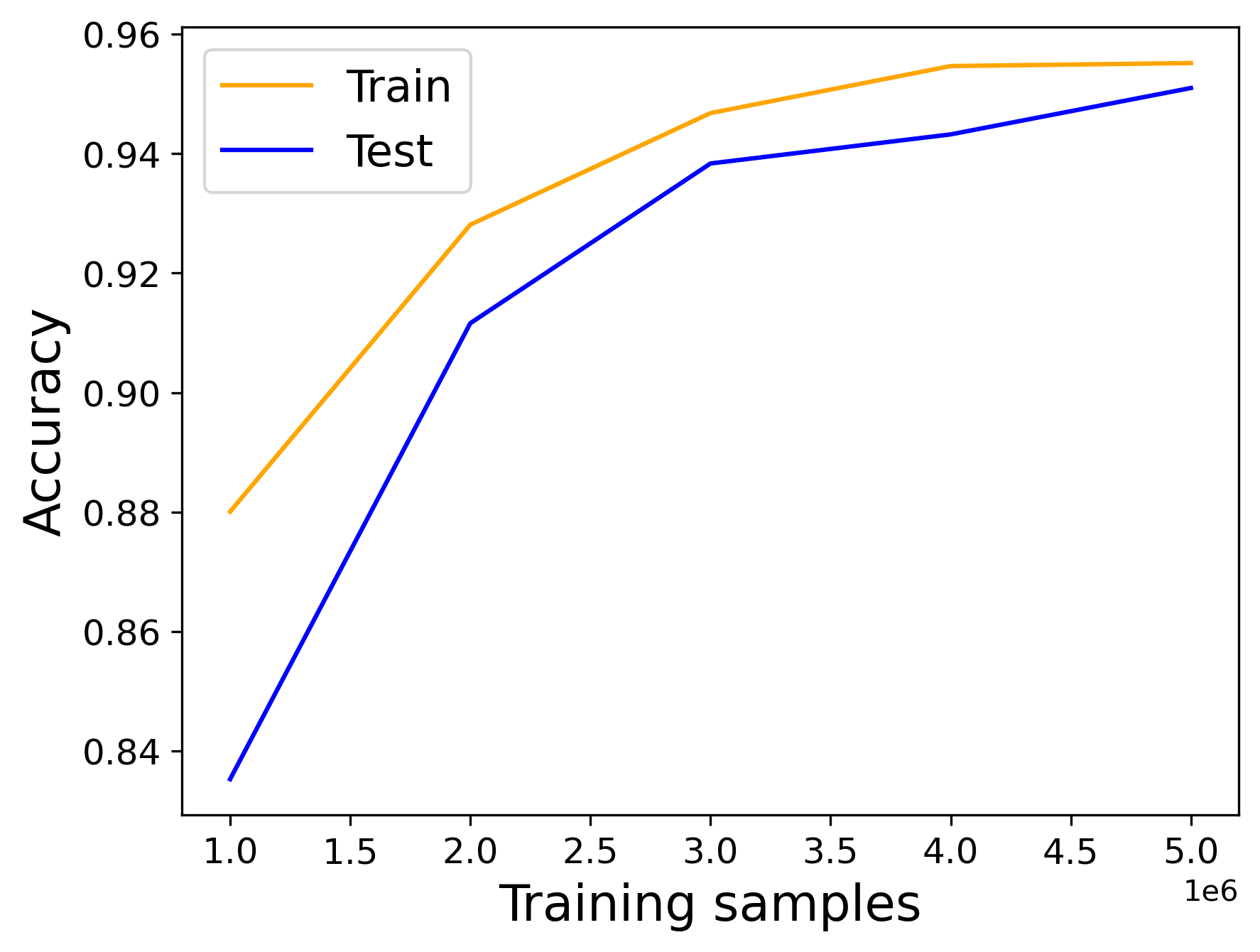}
        \caption{}\label{fig:learning_curve_train_size}
    \end{subfigure}
    \begin{subfigure}[t]{0.4\textwidth}
        \centering
        \includegraphics[width=\textwidth]{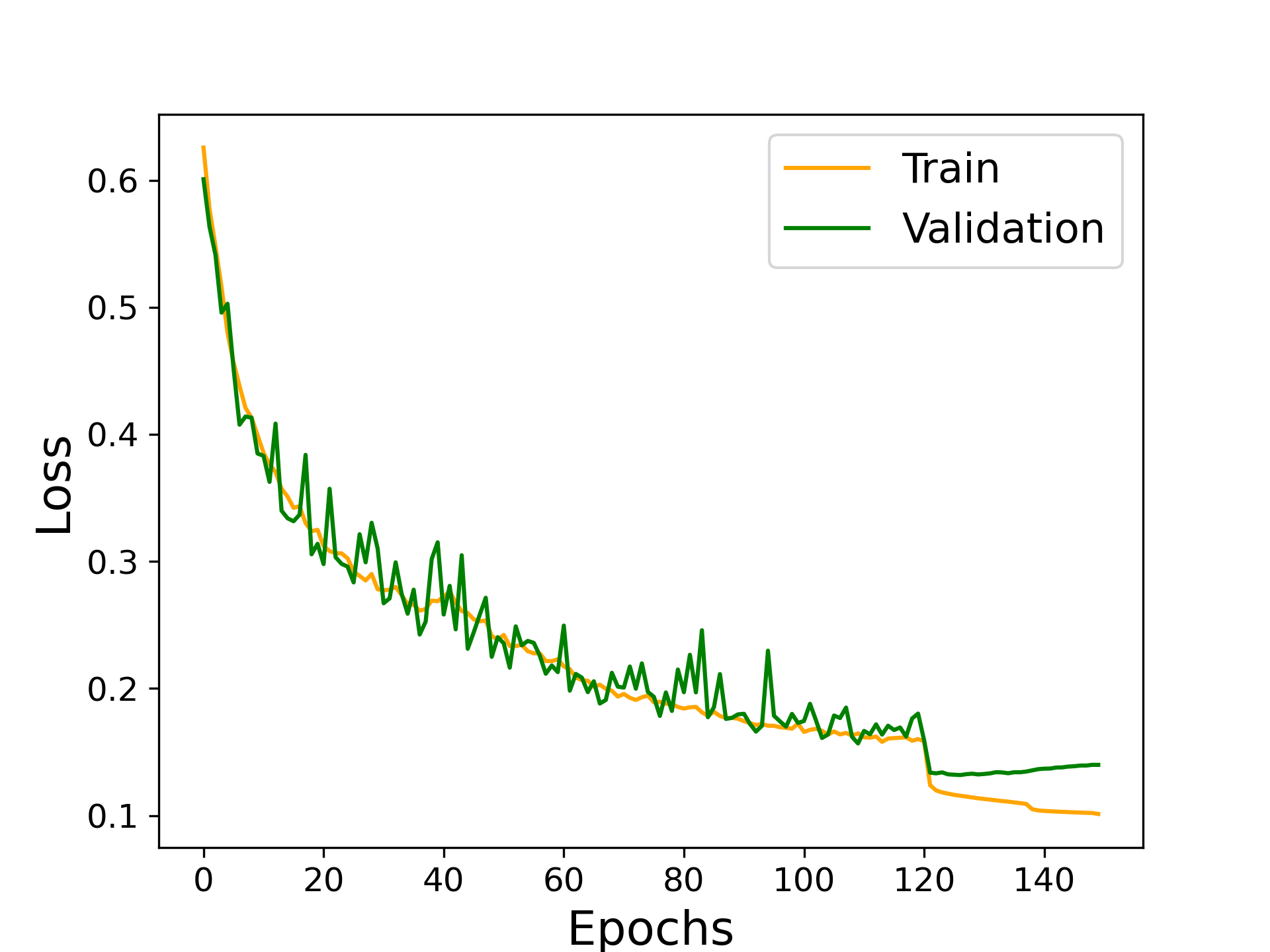}
        \caption{}\label{fig:learning_curve_terminality}
    \end{subfigure}
       \caption{(a) Accuracy for different train-test splits; (b) epochs against loss for the network trained on 5M~samples.}
\end{figure}

By trying different train-test splits, and using~20\% of the training data for validation throughout, we obtained the learning curve in Figure~\ref{fig:learning_curve_train_size}. This shows that a train-validate-test split of 4M-1M-5M produced an accurate model that did not overfit. Training this model gave the loss learning curve in Figure~\ref{fig:learning_curve_terminality}, and a final accuracy (on the test split of size~5M) of~$95\%$.

\section{Theoretical result}\label{sec:algorithm}
The high accuracy of the model in~\S\ref{sec:apply_ml} was very surprising. As explained in the introduction, \QFano{} varieties are of fundamental importance in algebraic geometry. However, asking whether a Fano variety has terminal singularities is, in general, an extremely challenging geometric question. In the case of a Fano toric variety one would typically proceed by constructing the fan, and then performing a cone-by-cone analysis of the combinatorics. This is computationally expensive and unsatisfying from a theoretical viewpoint. The success of the model suggested that a more direct characterisation is possible from the weight matrix alone. An analogous characterisation exists in the simpler case of weighted projective spaces~\cite{Kasprzyk2013}, which have Picard rank one, however no such result in higher Picard rank was known prior to training this model.

Inspired by this we prove a theoretical result, Proposition~\ref{general_result}, which leads to a new algorithm for checking terminality directly from the weight matrix, for \Qfactorial{} Fano toric varieties of Picard rank two. Consider a weight matrix as in~\eqref{eq:weight_matrix} that satisfies conditions~\eqref{item:zeros}--\eqref{item:well-formed} from~\S\ref{sec:data_generation}, and the toric variety~$X$ that it determines. As discussed in~\S\ref{sec:mathematical_background_main}, and explained in detail in~\S\ref{sec:mathematical_background_supp},~$X$ determines a convex polytope~$P$ in~$\RR^{N-2}$, with~$N$ vertices given by the first lattice points on the~$N$ rays of the fan. Each of the vertices of~$P$ is a lattice point (i.e., lies in~$\ZZ^{N-2} \subset \RR^{N-2}$), and~$X$ has terminal singularities if and only if the only lattice points in~$P$ are the vertices~$e_1,\ldots,e_N$ and the origin.

\begin{Defn}
    Let~$\Delta_i$ denote the simplex in~$\RR^{N-2}$ with vertices~$e_1, \dots, \hat{e}_i, \dots, e_{N}$ where~$e_i$ is omitted. We say that~$\Delta_i$ is~\emph{mostly empty} if each lattice point in~$\Delta_i$ is either a vertex or the origin.
\end{Defn}

\begin{Not}
    Let~$\{x\}$ denote the fractional part~$x - \lfloor x \rfloor$ of a rational number~$x$. 
\end{Not}

\begin{Prop}\label{general_result}
    Consider a weight matrix 
    \[
        \begin{bmatrix}
            a_1 & \cdots & a_{N} \\
            b_1 & \cdots & b_{N} 
        \end{bmatrix}
    \]         
    that satisfies conditions~\eqref{item:zeros}--\eqref{item:well-formed} from~\S\ref{sec:data_generation}. Let~$g_i = \gcd\{a_i, b_i\}$, and let~$A_i$,~$B_i$ be integers such that~$A_ia_i + B_ib_i =g_i$. Set
    \begin{align*}
        \alpha_i^j &= \frac{a_j b_i - b_j a_i}{g_i} & \alpha_i &= \sum_{j=1}^N \alpha_i^j \\
        \beta_i^j &= -A_ia_j - B_ib_j & \beta_i &= \sum_{j=1}^N \beta_i^j &
         f_i & = \frac{\alpha_i g_i}{\gcd\{g_i,\beta_i\}}
    \end{align*}
    noting that all these quantities are integers. Then~$\Delta_i$ is mostly empty if and only if for all~$k \in \{0, \dots, f_i-1\}$ and~$l \in \{0, \dots, g_i-1\}$ such that 
    \begin{align*}
        \sum_{j=1}^N \left\{ k\frac{\alpha_i^j}{f_i} + l\frac{\beta_i^j}{g_i}\right\} = 1
    \end{align*}
    we have that
    \[
        \left\{k \frac{\alpha_i^j}{f_i} + l\frac{\beta_i^j}{g_i}\right\} = \left\{ \frac{\alpha_i^j}{\alpha_i}\right\}
    \]
    for all~$j$.
\end{Prop}

Let~$s_+ = \{i \mid a_i b - b_i a > 0\}$, $s_- = \{i \mid a_i b - b_i a < 0 \}$, and let~$I$ be either~$s_+$ or~$s_-$. Then~$\Delta_i$, $i \in I$, forms a triangulation of~$P$. Thus~$X$ has terminal singularities if and only if~$\Delta_i$ is mostly empty for each~$i \in I$. This leads to Algorithm~\ref{alg:terminality_algorithm}.

\begin{algorithm} 
\caption{Test terminality for weight matrix $W=[[a_1, \dots, a_N],[b_1, \dots, b_N]]$.}\label{alg:terminality_algorithm}
\begin{algorithmic}[1]
    \State Set $a=\sum_{i=1}^Na_i$, $b=\sum_{i=1}^Nb_i$.
    \State Set $s_+ = \{i \mid a_i b - b_i a > 0\}$ and $s_- = \{i \mid a_i b - b_i a < 0 \}$.
    \State Set $I$ to be the smaller of $s_+$ and $s_-$.
\For{$i \in I$}
    \State Test if $\Delta_i$ is mostly empty, using Proposition~\ref{general_result}.
    \If{$\Delta_i$ is not mostly empty} 
        \State return \texttt{False}.
    \EndIf
\EndFor
\State return \texttt{True}.
\end{algorithmic}
\end{algorithm}

\subsection*{Comparisons}
Testing on 100\,000 randomly-chosen examples indicates that Algorithm~\ref{alg:terminality_algorithm} is approximately 15~times faster than the fan-based approach to checking terminality that we used when labelling our dataset (0.020s per weight matrix for Algorithm~\ref{alg:terminality_algorithm} versus 0.305s for the standard approach implemented in Magma). On single examples, the neural network classifier is approximately 30~times faster than Algorithm~\ref{alg:terminality_algorithm}. The neural network also benefits greatly from batching, whereas the other two algorithms do not: for batches of size~10\,000, the neural network is roughly 2000~times faster than Algorithm~\ref{alg:terminality_algorithm}.

\section{The terminal toric Fano landscape}\label{sec:application}
Having trained the terminality classifier, we used it to explore the landscape of \QFano{} toric varieties with Picard rank two. To do so, we built a large dataset of examples and analysed their~\emph{regularized quantum period}, a numerical invariant of \QFano{} varieties~\cite{Coatesetal2013}. For smooth low\nobreakdash-dimensional Fano varieties, it is known that the regularized quantum period is a complete invariant~\cite{CoatesCortiGalkinKasprzyk2016}. This is believed to be true in higher dimension, but is still conjectural. Given a \QFano{} variety~$X$, its regularized quantum period is a power series
\[
  \hat{G}_X (t) = \sum_{d=0}^{\infty} c_d t^d
\]
where~$c_0 = 1$, $c_1 = 0$, $c_d = d!\, r_d$, and~$r_d$ is the number of degree-$d$ rational curves in~$X$ that satisfy certain geometric conditions. Formally speaking,~$r_d$ is a degree-$d$, genus-zero Gromov--Witten invariant~\cite{KontsevichManin1997}. The~\emph{period sequence} of~$X$ is the sequence~$(c_d)$ of coefficients of the regularized quantum period. This sequence grows rapidly. In the case where~$X$ is a \QFano{} toric variety of Picard rank two, rigorous asymptotics for this growth are known. 

\begin{Thm}[Theorem 5.2,~\cite{CoatesKasprzykVeneziale2022}]
    Consider a weight matrix
    \[
        \begin{bmatrix}
            a_1 & \dots & a_N \\
            b_1 & \dots & b_N
        \end{bmatrix} 
    \]
    for a \Qfactorial{} Fano toric variety~$X$ of Picard rank two. Let~$a = \sum_{i=1}^N a_i$ and~$b = \sum_{i=1}^N b_i$, and let~$[\mu\!:\!\nu] \in \PP^1$ be the unique real root of the homogeneous polynomial 
    \begin{align}\label{eq:probability_equation}
        \prod_{i=1}^N (a_i \mu + b_i \nu)^{a_i b} -
        \prod_{i=1}^N (a_i \mu + b_i \nu)^{b_i a}
    \end{align}
    such that~$a_i \mu + b_i \nu \geq 0$ for all~$i \in \{1,2,\ldots,N\}$. Let~$(c_d)$ be the corresponding period sequence. Then non-zero coefficients~$c_d$ satisfy
    \[
        \log c_d \sim Ad - \frac{\dim{X}}{2} \log d + B
    \]
    as~$d \to \infty$, where 
    \begin{align}\label{eq:A_and_B}
        \begin{split}
            A &= -\sum_{i=1}^N p_i \log p_i \\
            B &= - \frac{\dim{X}}{2} \log(2 \pi) - \frac{1}{2} \sum_{i=1}^N \log p_i - \frac{1}{2} \log \left( \sum_{i=1}^N \frac{(a_i b - b_i a)^2}{ \ell^2 p_i} \right) 
        \end{split}
    \end{align}
    Here~$p_i = \displaystyle \frac{\mu a_i + \nu b_i}{\mu a + \nu b}$, so that~$\sum_i p_i = 1$, and~$\ell = \gcd\{a,b\}$ is the Fano index. 
\end{Thm}

In Figure~\ref{fig:PR2_hedgehog} we picture our dataset of \QFano{} varieties by using the coefficients~$A$ and~$B$ to project it to~$\RR^2$; for the corresponding images for terminal Fano weighted projective spaces, see~\cite[Figure~7a]{CoatesKasprzykVeneziale2022}. Note the stratification by Fano index. Although many weight matrices can give rise to the same toric variety, in our context we are using well-formed weight matrices in standard form~\eqref{eq:standard_form} and so at most two weight matrices can give rise to the same toric variety. We removed any such duplicates from our dataset, so the heatmap in Figure~\ref{fig:PR2_hedgehog_b} reflects genuine variation in the distribution of \QFano{} varieties, rather than simply the many-to-one correspondence between weight matrices and toric varieties.

\subsection*{Data generation}
The dataset pictured in Figure~\ref{fig:PR2_hedgehog} was generated using an AI\nobreakdash-assisted data generation workflow that combines algorithmic checks and our machine learning model, as follows.

\begin{itemize}
    \item Generate a random~$2\times10$ matrix with entries chosen uniformly from~$\{0,1,2,3,4,5,6,7\}$.
    \item Cyclically order the columns and only keep the matrix if it is in standard form, as in~\eqref{eq:standard_form}.
    \item Check conditions~\eqref{item:zeros}--\eqref{item:well-formed} from~\S\ref{sec:data_generation}.
    \item Predict terminality using the neural network classifier from~\S\ref{sec:apply_ml}, only keeping examples that are classified as terminal and storing their probabilities.
    \item Set~$\mu=1$ in~\eqref{eq:probability_equation} and solve the univariate real polynomial in the correct domain to obtain the solution~$(1,\nu)$.
    \item Calculate the coefficients~$A$ and~$B$ using the formulae in~\eqref{eq:A_and_B}.
\end{itemize}

The final dataset is composed of 100M samples. Each of these represents a \Qfactorial{} toric Fano variety of dimension eight and Picard rank two that the classifier predicts is a \QFano{} variety. 

\subsection*{Data analysis}
\begin{figure}
    \centering
    \begin{subfigure}[b]{0.47\textwidth}
        \centering
        \includegraphics[width=\textwidth]{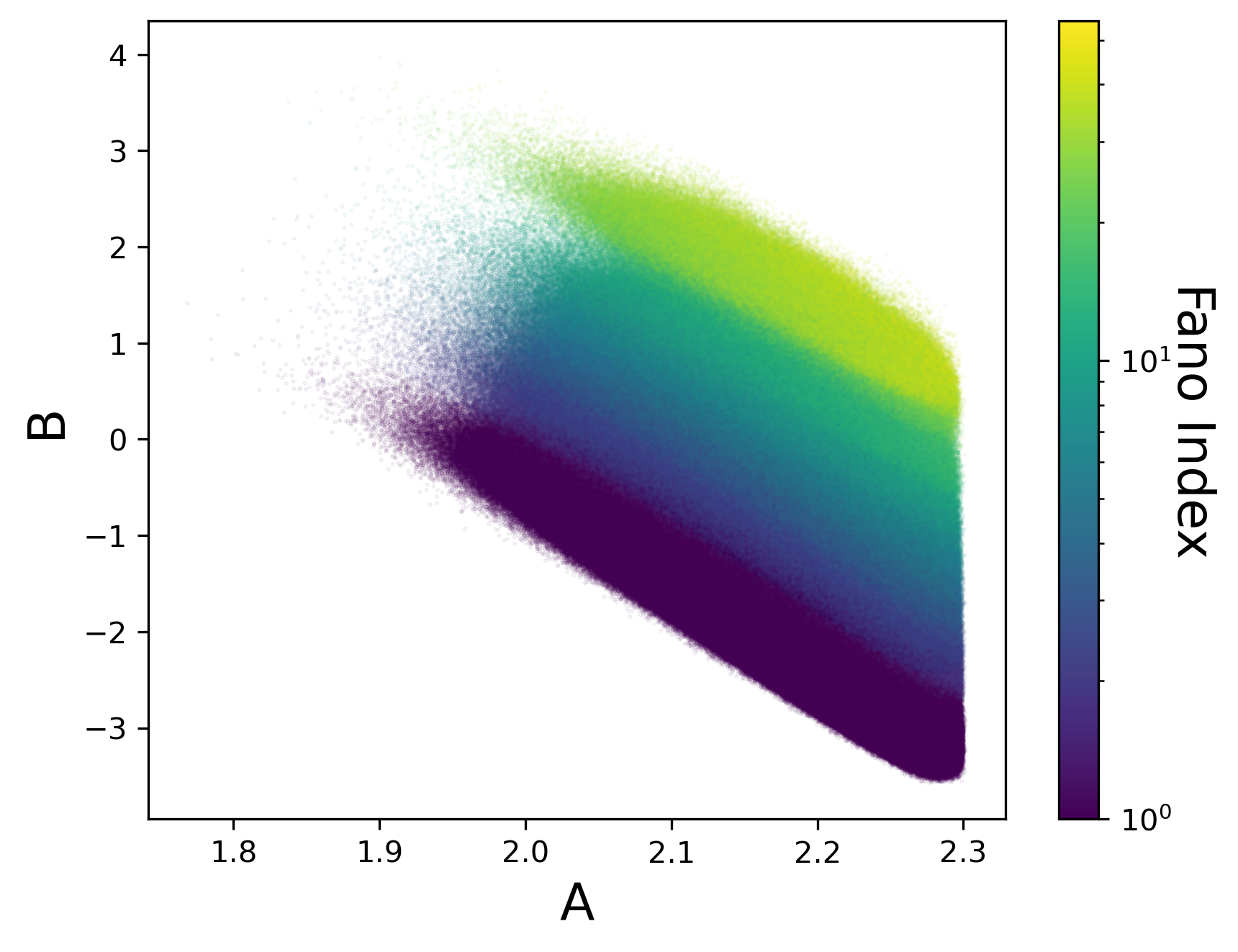}
        \caption{}\label{fig:PR2_hedgehog_a}
    \end{subfigure}
    \hfill
    \begin{subfigure}[b]{0.47\textwidth}
        \centering
        \includegraphics[width=\textwidth]{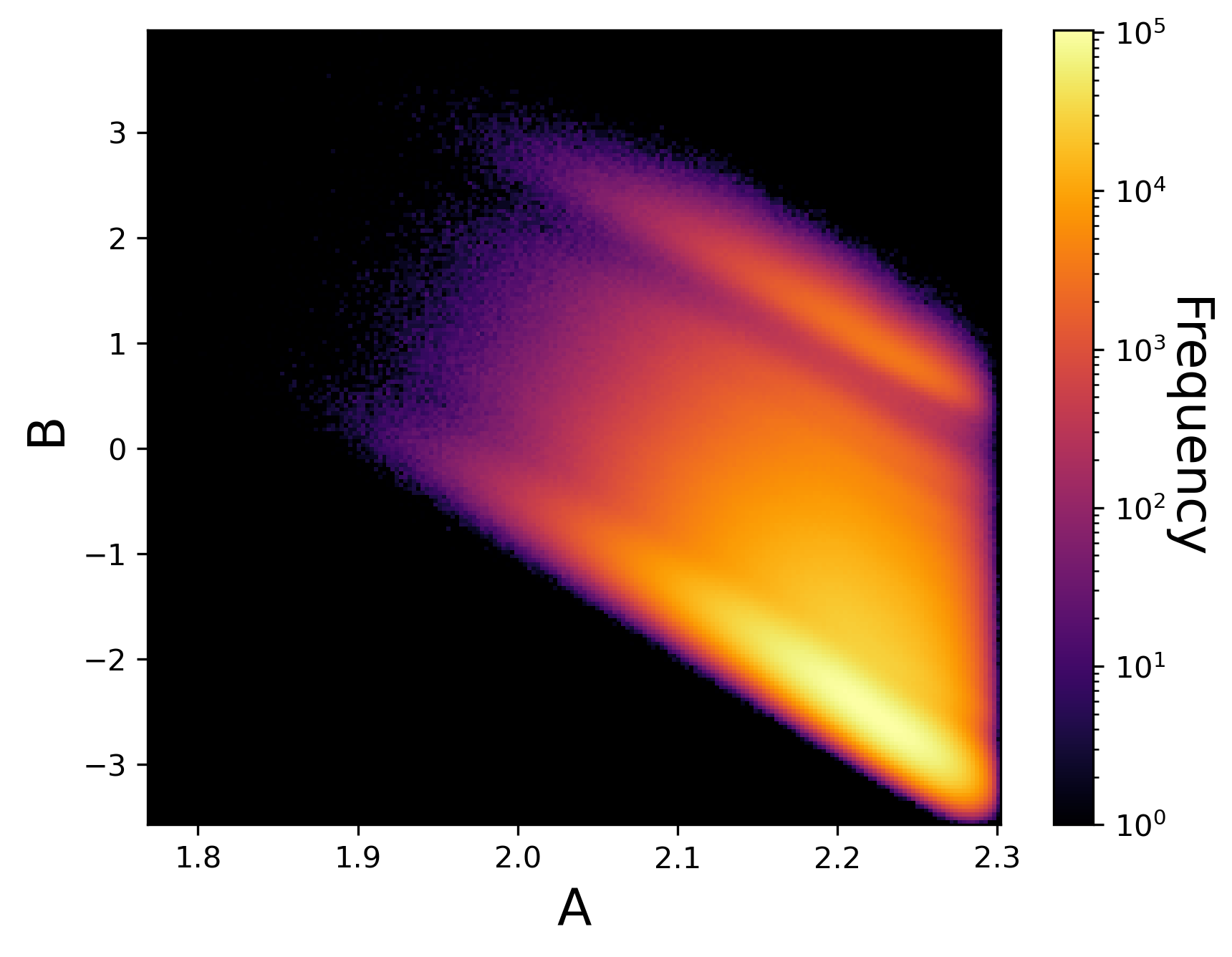}
        \caption{}\label{fig:PR2_hedgehog_b}
    \end{subfigure}
       \caption{A dataset of 100M probably-\QFano{} toric varieties of Picard rank two and dimension eight, projected to~$\RR^2$ using the growth coefficients~$A$ and~$B$ from~\eqref{eq:A_and_B}. In (a) we colour by Fano index, while in (b) we colour a heatmap according to the frequency.}
       \label{fig:PR2_hedgehog}
\end{figure}

We note that the vertical boundary in Figure~\ref{fig:PR2_hedgehog} is not a surprise. In fact, we can apply the log-sum inequality to the formula for~$A$ to obtain
\[
    A = -\sum_{i=1}^{N} p_i \log(p_i) \leq - \left(\sum_{i=1}^{N} p_i\right) \log \left(\frac{\sum_{i=1}^N p_i}{N} \right) = \log(N)
\]
In our case~$N = 10$, and the vertical boundary that we see in Figure~\ref{fig:PR2_hedgehog_a} is the line~$x = \log(10) \sim 2.3$. We also see what looks like a linear lower bound for the cluster; a similar bound was observed, and established rigorously, for weighted projective spaces in~\cite{CoatesKasprzykVeneziale2022}.

Closer analysis (see~\S\ref{sec:futher_data_analysis}) reveals large overlapping clusters that correspond to Fano varieties of different Fano index. Furthermore the simplest toric varieties of Picard rank two -- products of projective spaces, and products of weighted projective spaces -- appear to lie in specific regions of the diagram.

\section{Limitations and future directions}\label{sec:limitations}
\begin{figure}
    \centering
    \includegraphics[width=\textwidth]{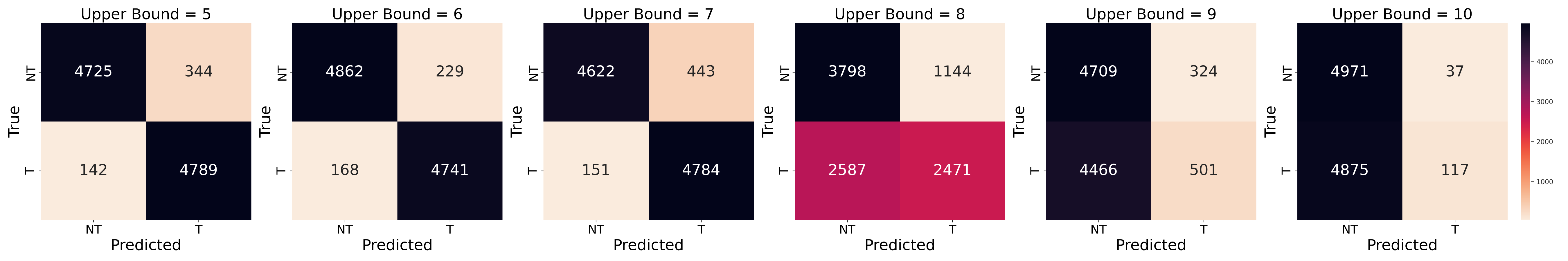}
    \caption{Confusion matrices for the neural network classifier on in-sample and out-of-sample data. In each case a balanced set of 10\,000 random examples was tested.}
       \label{fig:combined_confusion}
\end{figure}

The main message of this work is a new proposed AI\nobreakdash-assisted workflow for data generation in pure mathematics. This allowed us to construct, for the first time, an approximate landscape of objects of mathematical interest (\QFano{} varieties) which is inaccessible by traditional methods. We hope that this methodology will have broad application, especially to other large-scale classification questions in mathematics, of which there are many~\cite{LieAtlas,Cremona2016,Murmurations}.

Our approach has some limitations, however, which we enumerate here. Some of these limitations suggest directions for future research. A key drawback, common to most~ML models, is that our classifier performs poorly on out-of-sample data. Recall from~\S\ref{sec:data_generation} that the dataset we generated bounded the entries of the matrices by seven. For weight matrices within this range the model is extremely accurate (95\%), however this accuracy drops off rapidly for weight matrices that fall outside of this range: 62\% for entries bounded by eight; 52\% for entries bounded by nine; and 50\% for entries bounded by ten. See Figure~\ref{fig:combined_confusion} for details. Note that the network quickly degenerates to always predicting non-terminal singularities.

Furthermore the training process seems to require more data than we would like, given how computationally expensive the training data is to generate. It is possible that a more sophisticated network architecture, that is better adapted to this specific problem, might require less data to train.

Mathematically, our work here was limited to toric varieties, and furthermore only to toric varieties of Picard rank two. Finding a meaningful vectorisation of an arbitrary algebraic variety looks like an impossible task. But if one is interested in the classification of algebraic varieties up to deformation, this might be less of a problem than it first appears. Any smooth Fano variety in low dimensions is, up to deformation, either a toric variety, a toric complete intersection, or a quiver flag zero locus~\cite{CoatesCortiGalkinKasprzyk2016,Kalashnikov2019}; one might hope that this also covers a substantial fraction of the \QFano{} landscape. Each of these classes of geometry is controlled by combinatorial structures, and it is possible to imagine a generalisation of our vectorisation by weight matrices to this broader context.

Generalising to \Qfactorial{} Fano toric varieties in higher Picard rank will require a more sophisticated approach to equivariant machine learning. In this paper, we could rely on the fact that there is a normal form~\eqref{eq:standard_form} for rank-two weight matrices that gives an almost unique representative of each $\SL_2(\ZZ) \times S_N$\nobreakdash-orbit of weight matrices. For higher Picard rank~$r$ we need to consider weight matrices up to the action of~$G = \SL_r(\ZZ) \times S_N$. Here no normal form is known, so to work $G$\nobreakdash-equivariantly we will need to augment our dataset, to fill out the different $G$\nobreakdash-orbits, or to use invariant functions of the weights as features. The latter option, geometrically speaking, is working directly with the quotient space.

The best possible path forward would be to train an explainable model that predicted terminality from the weight data. This would allow us to extract from the machine learning not only that the problem is tractable, but also a precise mathematical conjecture for the solution. At the moment, however, we are very far from this. The multilayer perceptron that we trained is a black-box model, and post-hoc explanatory methods such as SHAP analysis~\cite{LundbergLee2017} yielded little insight: all features were used uniformly, as might be expected. We hope to return to this point elsewhere.

\subsection*{Data and code availability}
The datasets underlying this work and the code used to generate them are available from Zenodo under a CC0 license~\cite{dim8Data}. Data generation and post-processing was carried out using the computational algebra system Magma V2.27-3~\cite{BosmaCannonPlayoust1997}. The machine learning model was built using PyTorch v1.13.1~\cite{Pytorch} and scikit-learn v1.1.3~\cite{scikit-learn}. All code used and trained models are available from BitBucket under an MIT licence~\cite{supporting-code}.

\subsection*{Acknowledgements}
TC was partially supported by ERC Consolidator Grant~682603 and EPSRC Programme Grant~EP/N03189X/1. AK was supported by EPSRC Fellowship~EP/N022513/1. SV was supported by the Engineering and Physical Sciences Research Council [EP/S021590/1], the EPSRC Centre for Doctoral Training in Geometry and Number Theory (The London School of Geometry and Number Theory), University College London. The authors would like to thank Hamid Abban, Alessio Corti, and Challenger Mishra for many useful conversations, and the anonymous referees for their insightful feedback and suggestions.

\appendix
\section{Mathematical background}\label{sec:mathematical_background_supp}
\subsection*{Toric varieties}
The prototypical example of a toric Fano variety is two\nobreakdash-dimensional projective space, $\PP^2$. As mentioned in~\S\ref{sec:mathematical_background_main}, this is defined by taking the quotient of~$\CC^3\setminus \{\mathbf{0}\}$ by the following action of~$\CC^\times$:
\[
    \lambda \cdot (z_1, z_2, z_3) = (\lambda z_1, \lambda z_2, \lambda z_3 )
\]

The elements of~$\PP^2$ are equivalence classes that can be written as~$[z_1\!:\!z_2\!:\!z_3]$ where at least one of the~$z_i$ is non-zero. The algebraic variety~$\PP^2$ is~\emph{smooth}, since we can cover it by three open subsets that are each isomorphic to the complex plane~$\CC^2$. Namely, 
\begin{align*}
    U_1 &= \{[z_1\!:\!z_2\!:\!z_3] \in \PP^2 \mid z_1 \neq 0\} \\
    U_2 &= \{[z_1\!:\!z_2\!:\!z_3] \in \PP^2 \mid z_2 \neq 0\}\\
    U_3 &= \{[z_1\!:\!z_2\!:\!z_3] \in \PP^2 \mid z_3 \neq 0\}
\end{align*}
To see that~$U_1$ is isomorphic to~$\CC^2$, we note that since~$z_1 \neq 0$ it can be rescaled to one. Therefore, each point in~$U_1$ can be identified with a (unique) point of the form~$[1\!:\!\bar{z}_2\!:\!\bar{z}_3]$; this gives the isomorphism to~$\CC^2$. Similar arguments show that~$U_2$ and~$U_3$ are each isomorphic to~$\CC^2$.

More generally, $(N-1)$\nobreakdash-dimensional projective space~$\PP^{N-1}$ is smooth, since it can be covered by~$N$ open subsets each isomorphic to~$\CC^{N-1}$. By modifying the action of~$\CC^{\times}$ on~$\CC^{N} \setminus \{\mathbf{0}\}$ we can define more general examples of toric varieties,~\emph{weighted projective spaces}, which in general contain singular points.

For example, we can consider the action of~$\CC^{\times}$ on~$\CC^3 \setminus \{\mathbf{0}\}$ defined by
\[
    \lambda \cdot (z_1, z_2, z_3) = (\lambda z_1, \lambda z_2, \lambda^2 z_3 )
\]
which gives rise to the weighted projective space~$\PP(1,1,2)$. Here the entries of the vector~$(1,1,2)$ are called the~\emph{weights} of the variety. In order to see that this variety is not smooth, we can consider the same open sets as above, 
\begin{align*}
    U_1 &= \{[z_1\!:\!z_2\!:\!z_3] \in \PP^2 \mid z_1 \neq 0\} \\
    U_2 &= \{[z_1\!:\!z_2\!:\!z_3] \in \PP^2 \mid z_2 \neq 0\}\\
    U_3 &= \{[z_1\!:\!z_2\!:\!z_3] \in \PP^2 \mid z_3 \neq 0\}
\end{align*}
As before,~$U_1$ and~$U_2$ are each isomorphic to~$\CC^2$. However,~$U_3$ is not. In fact, since~$z_3 \neq 0$ we can rescale the last entry to one, but the square in the definition of the action implies that there are two ways of doing so:
\[
    \pm z_3^{-1/2} \cdot (z_1, z_2, z_3) = (\pm z_3^{-1/2} z_1, \pm z_3^{-1/2} z_2, 1) 
\] 
Therefore,~$U_3 \cong \CC^2/\mu_2$ where~$\mu_2 = \{1,-1\}$ is the group of square roots of unity. Note that~$\CC^2/\mu_2$ has a singular point at the origin, which corresponds to the singular point~$[0\!:\!0\!:\!1]$ in~$U_3$. We say that~$\PP(1,1,2)$ has two smooth charts,~$U_1$ and~$U_2$, and one singular chart~$U_3$.

This generalises to higher dimensions by considering~$\CC^{\times}$ acting on~$\CC^N \setminus \{\mathbf{0}\}$ by 
\[
    \lambda \cdot (z_1,\dots, z_N) = (\lambda^{a_1} z_1, \dots , \lambda^{a_N} z_N )
\]
for some choice of weights~$(a_1, \dots, a_N) \in \ZZ^{N}_{>0}$. The algebraic variety~$\PP(a_1,a_2,\ldots,a_N)$ is an $(N-1)$\nobreakdash-dimensional \Qfactorial{} Fano toric variety of Picard rank one, called a~\emph{weighted projective space}~\cite{IanoFletcher2000,Dolgachev1982}. Setting the~$a_i$ equal to~1 recovers~$\PP^{N-1}$.

For any two weighted projective spaces~$X=\PP(a_1,\ldots,a_N)$ and~$Y=\PP(b_1,\ldots,b_M)$, we can consider their product~$X \times Y$. This arises as a quotient of~$\CC^{N+M}$ by an action of~$\CC^{\times} \times \CC^{\times}$, where the first~$\CC^{\times}$ acts on the first~$N$ co-ordinates of~$\CC^{N+M}$ and the second~$\CC^{\times}$ acts on the last~$M$ co-ordinates. The two actions are specified by the weights of each weighted projective space. We can summarise this information in a~\emph{weight matrix}
\begin{align*}
    \begin{bmatrix}
        a_1 & \cdots & a_{N} & 0 & \cdots & 0  \\
        0 & \cdots & 0 & b_1 & \cdots & b_{M}
    \end{bmatrix} 
\end{align*}
This type of construction can be generalised to any action of~$\CC^{\times} \times \CC^{\times}$ on~$\CC^N$ given defined as 
\[
    (\lambda, \mu) \cdot (z_1, \dots, z_{N}) = (\lambda^{a_1} \mu^{b_1} z_1 , \dots, \lambda^{a_{N}} \mu^{b_{N}} z_{N} )
\] 
and which can be encoded in a weight matrix of the form
\begin{align*}
    \begin{bmatrix}
        a_1 & \cdots & a_{N} \\
        b_1 & \cdots & b_{N}
    \end{bmatrix} 
\end{align*}
Note that in the case of projective spaces and weighted projective spaces we were considering~$\CC^N \setminus \{\mathbf{0}\}$, excluding the origin because it lies in the closure of every orbit. When generalising this concept, we need to exclude more points than just the origin for the quotient to be reasonable; explicitly we consider~$\CC^N \setminus S$, where~$S= S_+ \cup S_-$ for linear subspaces
\begin{equation*}
    \begin{aligned}
        &S_+ = \{(z_1, \dots, z_{N}) \mid z_i =0 \text{ if } b_i/a_i > b/a\} \\
        &S_- = \{(z_1, \dots, z_{N}) \mid z_i =0 \text{ if } b_i/a_i < b/a\} 
    \end{aligned}
\end{equation*}
and~$a = \sum_{i=1}^{N} a_i$,~$b = \sum_{i=1}^{N} b_i$: see~\cite{BrownCortiZucconi2004}. The resulting quotient~$X= (\CC^N \setminus S)/(\CC^{\times})^2$ is an $(N-2)$\nobreakdash-dimensional toric variety. If the linear subspaces~$S_+$ and~$S_-$ each have dimension at least two then~$X$ has Picard rank two.

\subsection*{From weight matrices to fans}
In~\S\ref{sec:mathematical_background_main}, a toric variety~$X$ was determined by a matrix
\begin{align}\label{eq:weight_matrix_supp}
    W = \begin{bmatrix}
        a_1 & \cdots & a_{N} \\
        b_1 & \cdots & b_{N}
    \end{bmatrix} 
\end{align}
that, as recalled above, records the weights of an action of~$(\Cstar)^2$ on~$\CC^N$. We will now explain how to recover the fan~$\Sigma(X)$ for the toric variety from this data~\cite{Fulton1993,CoxLittleSchenk2011}.

Consider the right kernel of the matrix~$W$, regarded as a~$\ZZ$-linear map. The kernel is a free submodule of~$\ZZ^N$, of rank~$N-2$, and choosing a basis for this submodule defines an~$N \times (N-2)$ matrix~$M$ such that~$W M = 0$. The rows of~$M$ define distinct primitive vectors~$e_1,\ldots,e_N$ in~$\ZZ^{N-2}$ such that
\begin{equation*}
    \begin{split}
            a_1 e_1 +\cdots +a_{N} e_{N} &= 0   \\
            b_1e_1 +\cdots +b_{N}e_{N} &= 0
    \end{split}
\end{equation*}
By construction, the vectors~$e_1,\ldots, e_N$ span the kernel of~$W$ over~$\ZZ$.

In general the construction of a toric variety (or equivalently a fan) from a weight matrix depends also on the choice of a~\emph{stability condition}, which is an element~$\omega$ of the column space of~$W$. In our case, however, because~$X$ is Fano there is a canonical choice for~$\omega$ given by~$(a,b)$, the sum of the columns of~$W$.  Let us denote the~$i$th column of~$W$ by~$D_i$. We set
\[
    \mathcal{A}_{\omega} = \left\{I \subset\{1,2, \dots, N\} \mid \omega \in \angle_I\right\} 
\]
where
\[
    \angle_I = \left\{\sum_{i \in I} \lambda_i D_i\ \Big|\ \lambda_i \in \RR_{>0}\right\} 
\]
The fan~$\Sigma(X)$ is the collection of cones in~$\RR^{N-2}$ given by
\begin{align*}
    \{\sigma_I \mid \bar{I}  \in \mathcal{A}_\omega\} && \text{where $\sigma_I = \Cone\{e_i \mid i \in I\}$}
\end{align*}
Here~$\bar{I}$ is the complement of~$I$ in~$\{1,2,\ldots,N\}$. 

Recall our assumptions on the weight matrix~$W$:
\begin{enumerate}
    \setcounter{enumi}{-1}
    \item~\label{item:convex} The columns of~$W$ span a strictly convex cone in~$\RR^2$.
    \item~\label{item:zeros_supp} None of the columns are the zero vector.
    \item~\label{item:simplicial_supp} The sum of the columns is not a multiple of any of them.
    \item~\label{item:rays_supp} The subspaces~$S_+$ and~$S_-$, defined in~\eqref{eq:S_plus_minus}, are both of dimension at least two.
\end{enumerate}
(We number from zero here to match the numbering of the conditions in~\S\ref{sec:data_generation}.)
Conditions~\eqref{item:convex} and~\eqref{item:zeros_supp} together guarantee that the fan~$\Sigma(X)$ is complete; that is, its support covers~$\RR^{N-2}$. The toric variety~$X$ is therefore compact. Condition~\eqref{item:simplicial_supp} ensures that each top\nobreakdash-dimensional cone in the fan has~$N-2$ rays; that is, the fan is simplicial. This implies that the toric variety~$X$ is \Qfactorial{}. Condition~\eqref{item:rays_supp} ensures that each of the vectors~$e_1,\ldots,e_N$ generates a one\nobreakdash-dimensional cone~$\RR_{\geq 0} e_i$ in the fan~$\Sigma(X)$. Together with \Qfactorial{}ity, this implies that the Picard rank of~$X$ is two.

\subsection*{Checking terminality} 
Each top\nobreakdash-dimensional cone~$\sigma$ in~$\Sigma(X)$ is generated over~$\RR_{\geq 0}$ by~$N-2$ of the vectors~$e_1,\ldots,e_N$. These generators are contained in a unique $(N-3)$\nobreakdash-dimensional hyperplane~$H$. The cone~$\sigma$ corresponds to a terminal singularity in~$X$ if and only if the only lattice points in~$\sigma$ that lie on or below~$H$ are the generators of~$\sigma$ and the origin~\cite{Reid1987}. $X$ has terminal singularities if and only if each top\nobreakdash-dimensional cone of~$\Sigma(X)$ corresponds to a terminal singularity. This justifies the assertion, given in~\S\ref{sec:mathematical_background_main}, that~$X$ has terminal singularities if and only if the convex polytope~$P = \conv\{e_1,\ldots,e_N\}$ is mostly empty.

\subsection*{A subtlety with quotient gradings}
In~\S\ref{sec:introduction}, in the paragraph `Why dimension eight?', we noted that the analogue of our dataset in dimension three contains 34~examples. There are~35 \QFano{} toric varieties of Picard rank two in dimension three~\cite{Kasprzyk2006}, but precisely one of these has a quotient grading and so does not fit into the framework we consider here. The exception is~$X = \PP^1\times\PP^2 / \mu_3$, where~$\mu_3$ acts via~$(u,v;x,y,z)\mapsto (u,\varepsilon v;x,\varepsilon y,\varepsilon^2z)$ and~$\varepsilon$ is a primitive cube root of unity. The quotient grading arises here because the primitive generators for rays of the fan~$\Sigma(X)$ fail to span the ambient lattice over~$\ZZ$. If we instead regard the primitive generators as living inside the sublattice that they generate, then we recover one of the other 34~terminal examples:~$\PP^1 \times \PP^2$. The analogue of this phenomenon happens in higher dimensions too, and so we ignore quotient gradings in our methodology.

\subsection*{Significance of \QFano{} varieties}
As mentioned in~\S\ref{sec:introduction}, \QFano{} varieties are `atomic pieces' from which more complicated algebraic varieties are made, and so one can think of the classification of \QFano{} varieties as building a Periodic Table for geometry. Understanding this classification is a fundamental problem in algebraic geometry, and is the motivation behind a huge amount of research; see e.g.~\cite{Cascini2021,KollarMori1998,Kollar1987,GRDB} and the references therein. 

\QFano{} varieties also play an important role elsewhere in mathematics, for example in the study of K-stability and the existence of K\"ahler--Einstein metrics~\cite{Berman2016}. In theoretical physics, \QFano{} varieties provide, through their `anticanonical sections', the main construction of the Calabi-Yau manifolds which give geometric models of spacetime~\cite{Polchinski2005, Greene1997, CandelasHorowitzStromingerWitten1985} in Type II string theory.

Moreover, terminal singularities -- the focus of this paper -- are the singularities that appear in the Minimal Model Program~\cite{Kollar1987}, and they also occur across mathematics. For example, in F-theory, terminal singularities reflect the presence of localized matter states from wrapped M2-branes which are not charged under any massless gauge potential~\cite{ArrasGrassiWeigand2018}. Moreover, in the toric context, having only terminal singularities means that the corresponding polytope contains no lattice points other than the origin and the vertices. These are referred to in the combinatorics literature as one-point lattice polytopes, and are important in optimisation problems.

\section{Further data analysis}\label{sec:futher_data_analysis}
\begin{figure}[ht]
     \centering
     \begin{subfigure}{0.4\textwidth}
         \includegraphics[width=1.0\linewidth]{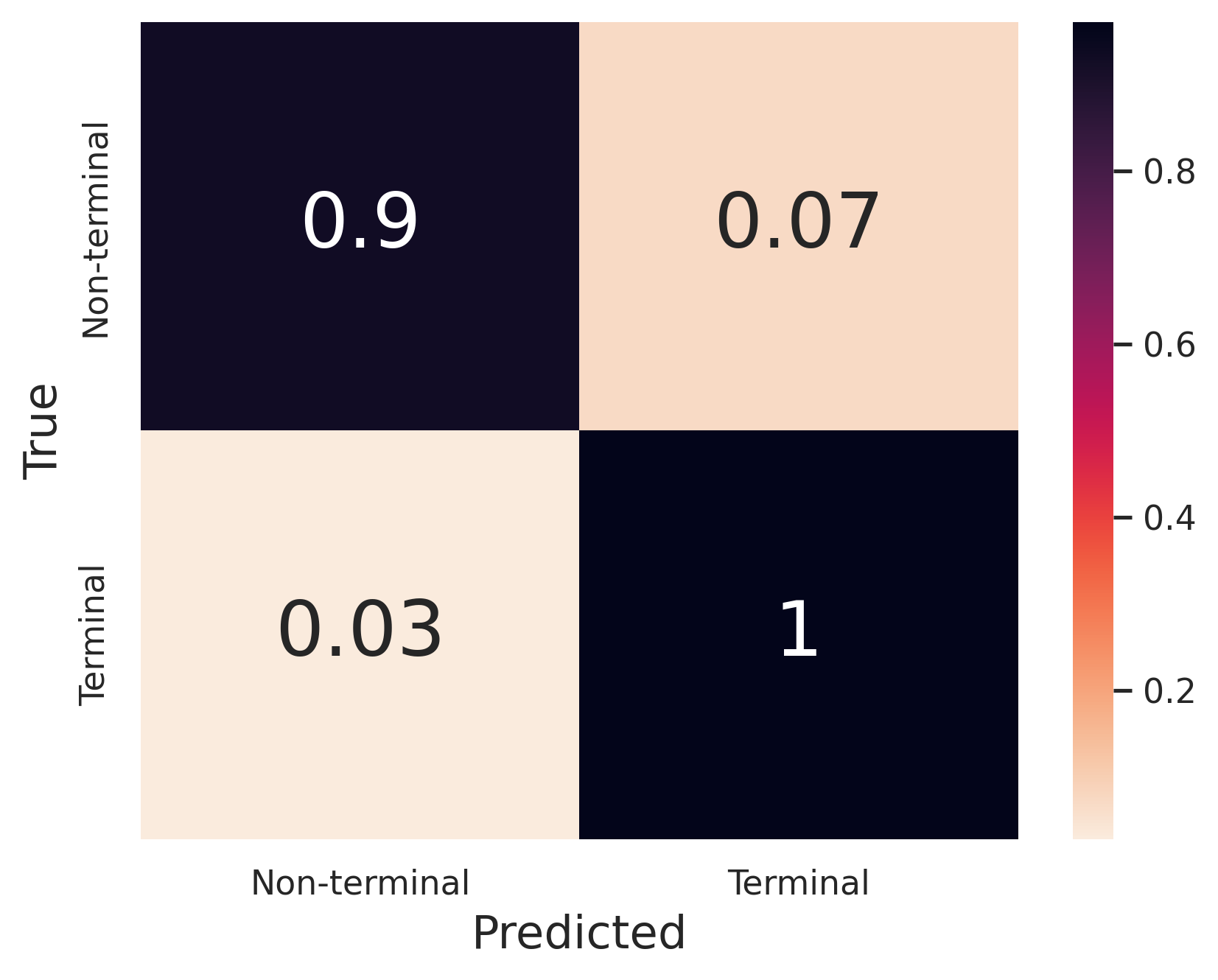}
         \caption{}
     \end{subfigure}
     \begin{subfigure}{0.4\textwidth}
         \includegraphics[width=1.0\linewidth]{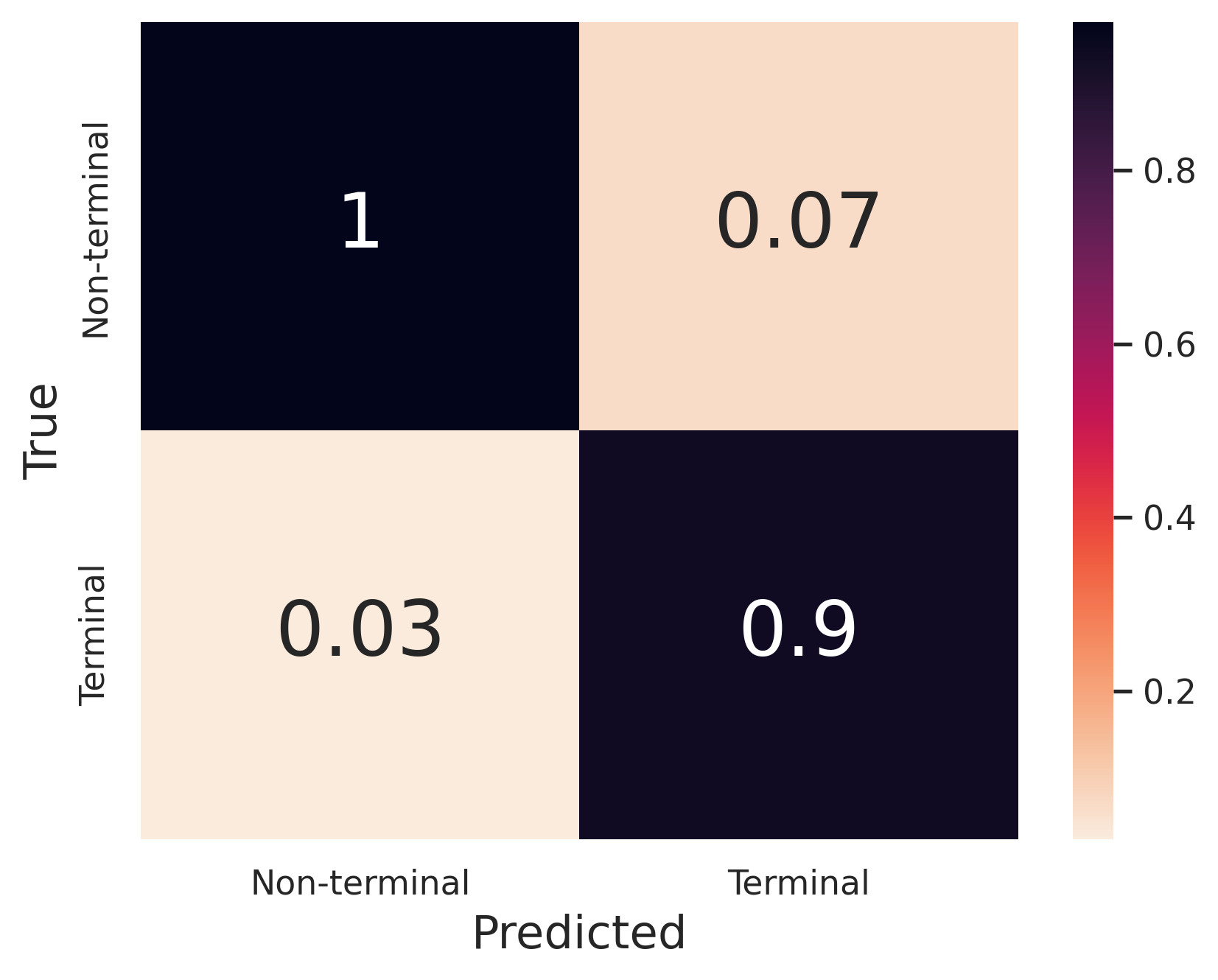}
         \caption{}
     \end{subfigure}
     \caption{Confusion matrices for the classifier trained on 5M~samples: (a)~is normalised with respect to the true axis; (b)~is normalised with respect to the predicted axis.}\label{fig:confusion_matrices_terminality}
 \end{figure}
 
The neural network classifier described in~\S\ref{sec:apply_ml} is remarkably accurate at determining whether a \Qfactorial{} Fano toric variety of Picard rank two and dimension eight is terminal or not. Confusion matrices for the classifier are presented in Figure~\ref{fig:confusion_matrices_terminality}. Because of this high accuracy, we were able to use this classifier to generate a dataset of 100M probably-\QFano{} toric varieties of Picard rank two and dimension eight; see~\S\ref{sec:application}. Creating this first glimpse of the \QFano{} landscape would have been impractical using conventional methods. Based on the timing data outlined in~\S\ref{sec:computational_resources} below, we estimate that generating this dataset using conventional methods would have taken 160~days on our HPC cluster, equivalent to 600~CPU~\emph{years}. In contrast, by using the neural network classifier and batch processing we were able to generate this dataset in under 120~CPU~\emph{hours}.

One striking feature of the landscape of 100M probably-\QFano{} toric varieties, plotted in Figure~\ref{fig:PR2_hedgehog}, is the stratification by Fano index. Recall that the Fano index of~$X$ is equal to the greatest common divisor of~$a$ and~$b$, where~$(a,b)$ is the sum of the columns of the matrix~\eqref{eq:weight_matrix_supp}. For our dataset, the entries in the matrix~\eqref{eq:weight_matrix_supp} are bounded between~zero and~seven, and hence the range of possible Fano indices that can appear in the dataset is bounded. Figure~\ref{fig:PR2_hedgehog} appears to show overlapping clusters of cases, with the Fano index increasing as we move from the bottom of the plot (Fano index one) to the top.

\subsection*{Products of weighted projective space}
To better understand this clustering by Fano index, we consider the simplest \Qfactorial{} Fano toric varieties of Picard rank two: products of weighted projective spaces. Recall from~\S\ref{sec:mathematical_background_supp} that a product of weighted projective spaces~$X=\PP(a_1,\ldots,a_N)$ and~$Y=\PP(b_1,\ldots,b_M)$ is specified by a weight matrix 
\[
\begin{bmatrix}
a_1&\cdots&a_N&0&\cdots&0\\
0&\cdots&0&b_1&\cdots&b_M
\end{bmatrix}
\]
This matrix determines a \Qfactorial{} Fano toric variety of Picard rank two and dimension~$N+M-2$, denoted~$X\times Y$. The singular points of~$X\times Y$ are determined by the singular points of~$X$ and~$Y$. In particular,~$X\times Y$ is terminal if and only if both~$X$ and~$Y$ are terminal.

In general a weighted projective space~$X=\PP(a_1,a_2,\ldots,a_N)$ may have singular points; these are determined by the weights~$(a_1,a_2,\ldots,a_N)$. Proposition~2.3 of~\cite{Kasprzyk2013} characterises when the singular points of~$X$ are terminal. Namely,~$X$ is terminal if and only if
\[
\sum_{i=1}^N\{ka_i/a\}\in\{2,\ldots,N-2\}
\]
for each~$k\in\{2,\ldots,a-2\}$. Here~$a=a_1+a_2+\cdots+a_N$, and~$\{x\}$ denotes the fractional part~$x-\lfloor x\rfloor$ of a rational number~$x$. This is the Picard rank one analogue to Proposition~\ref{general_result}.

\begin{figure}[ht]
    \centering
    \begin{subfigure}{0.45\textwidth}
        \includegraphics[width=1.0\linewidth]{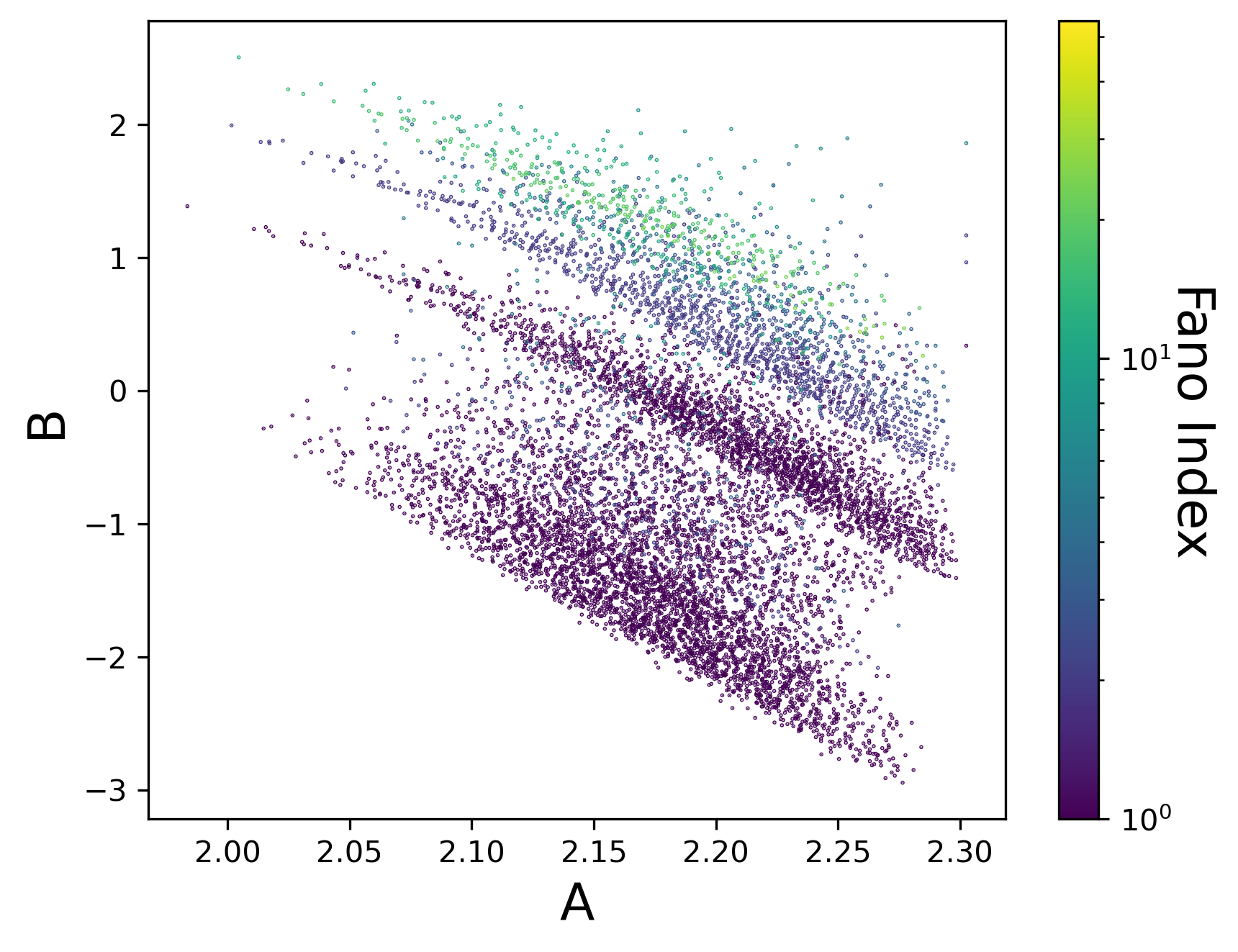}
        \caption{}
    \end{subfigure}
    \begin{subfigure}{0.45\textwidth}
        \includegraphics[width=1.0\linewidth]{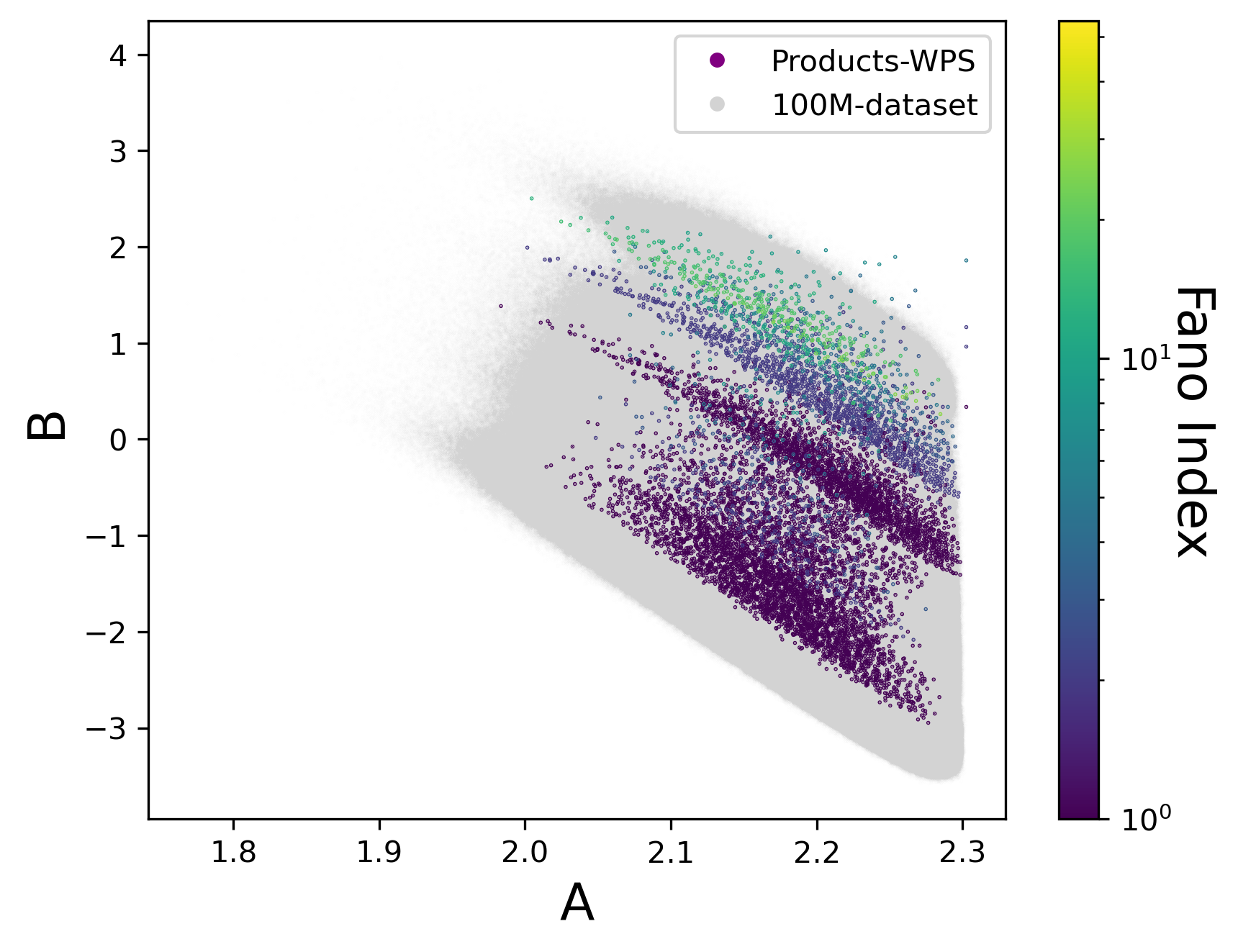}
        \caption{}
    \end{subfigure}
    \caption{\QFano{} products of weighted projective space in dimension eight, with weights bounded by~seven. (a)~Projection to $\RR^2$ using the growth coefficients from~\eqref{eq:A_and_B}. (b)~The same as~(a), but plotted on top of the dataset of 100M probably-\QFano{} toric varieties, plotted in grey.}\label{fig:prod_wps}
\end{figure}

We can enumerate all terminal weighted projective spaces in dimensions one to seven, with weights~$1\leq a_i\leq 7$, using the characterisation of terminal weighted projective space described above. The number in each dimension is given in Table~\ref{tab:terminal_wps}. By taking products, we obtain~8792 distinct \QFano{} toric varieties of Picard rank two in dimension eight; these examples are plotted in Figure~\ref{fig:prod_wps}. This supports our observation that the \QFano{} varieties fall into large overlapping clusters that are determined by the Fano index. Note that the products of weighted projective space appear to fall within the upper region of each cluster.

\begin{table}[t]
    \centering
    \vspace{1em}
    \begin{tabular}{r|ccccccc}
        \toprule
        $d$&1&2&3&4&5&6&7\\
        \#&1&1&7&80&356&972&2088\\
        \bottomrule
    \end{tabular}
    \caption{The number of terminal weighted projective spaces in dimension~$d$, $1\leq d\leq 7$, with weights $a_i$ bounded by~seven.}\label{tab:terminal_wps}
\end{table} 

\begin{figure}[ht]
    \centering
    \begin{subfigure}{0.45\textwidth}
        \includegraphics[width=1.0\linewidth]{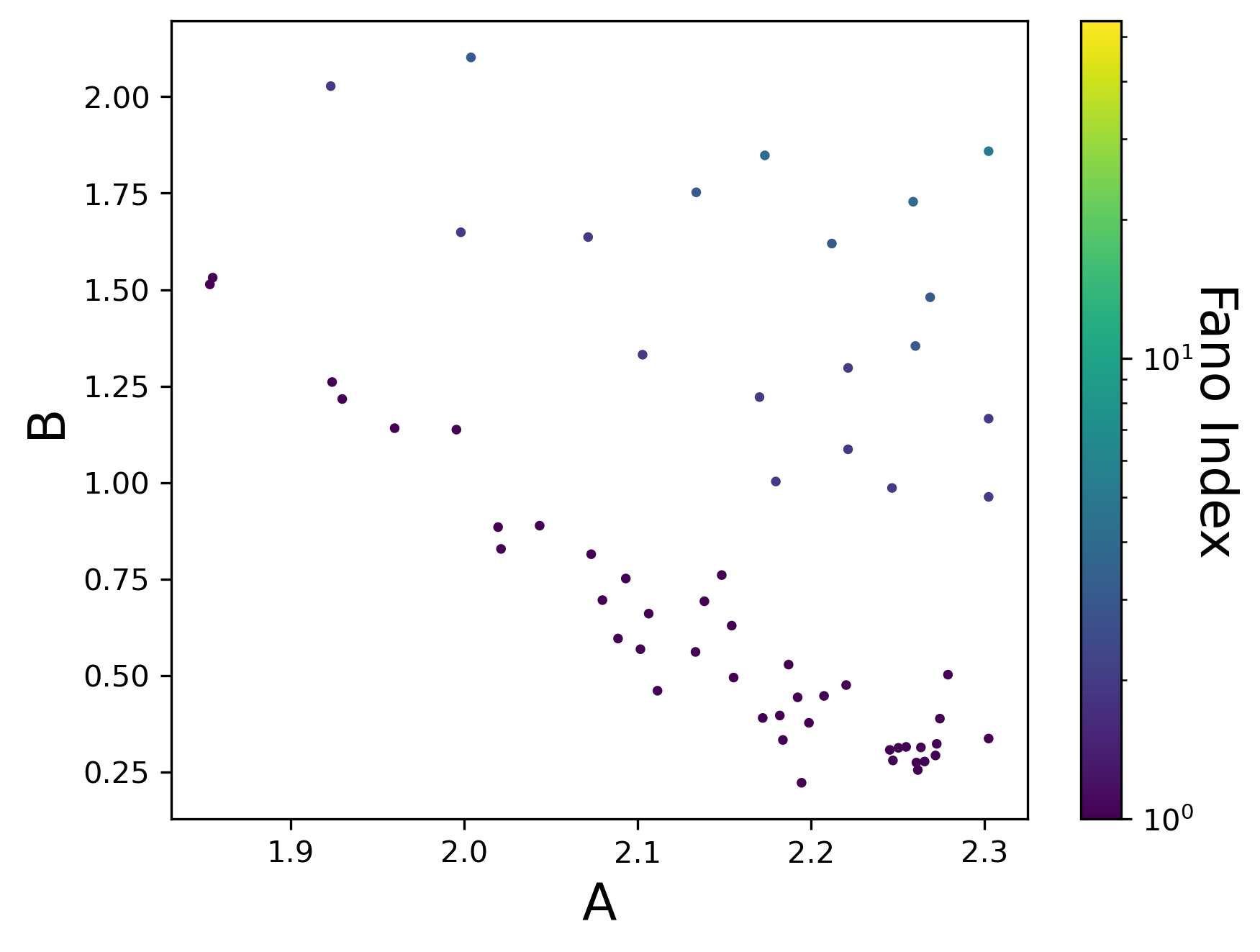}
        \caption{}
    \end{subfigure}
    \begin{subfigure}{0.45\textwidth}
        \includegraphics[width=1.0\linewidth]{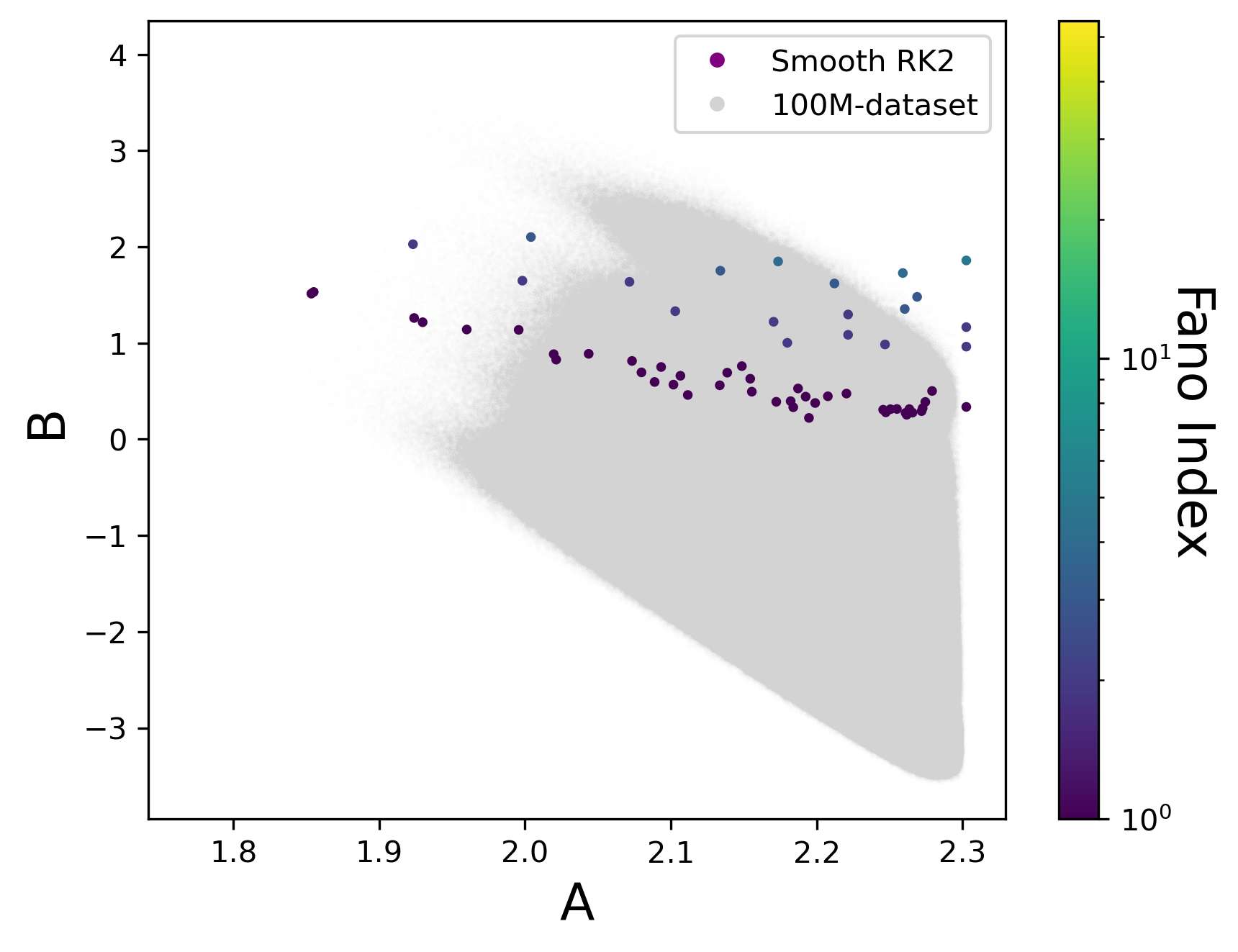}
        \caption{}
    \end{subfigure}
    \caption{The smooth Fano toric varieties in dimension eight and of Picard rank two. (a)~Projection to $\RR^2$ using the growth coefficients from~\eqref{eq:A_and_B}. (b)~The same as~(a), but plotted on top of the dataset of 100M probably-\QFano{} toric varieties, plotted in grey.}\label{fig:smooth}
\end{figure}

\subsection*{Smooth Fano toric varieties}
Projective space~$\PP^{N-1}$ is smooth, and so products of projective spaces are also smooth. More generally, the smooth Fano toric varieties up to dimension eight have been classified~\cite{Obro2007}. There are 62~smooth Fano toric varieties in dimension eight and of Picard rank two, all of which have weights bounded by~seven when expressed in standard form~\eqref{eq:standard_form}. These are plotted in Figure~\ref{fig:smooth}, and appear to fall in the upper extreme region within each cluster.

\begin{figure}[ht]
    \centering
    \includegraphics[width=0.45\textwidth]{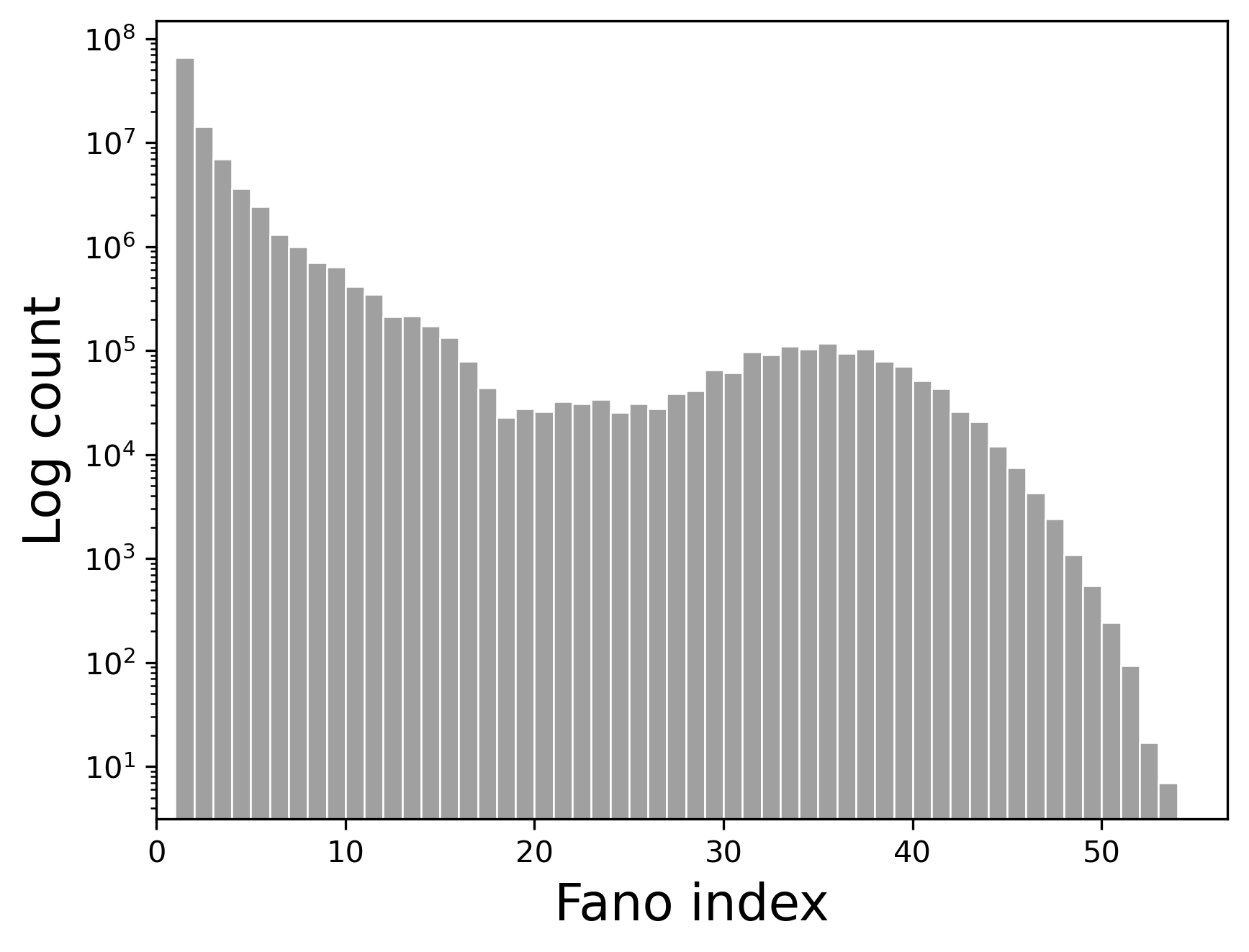}
    \caption{Distribution of the Fano index $\gcd\{a,b\}$ in the dataset of 100M probably-\QFano{} toric varieties (note that the vertical axis scale is logged).}
    \label{fig:indices_histogram}
\end{figure}

\begin{figure}[ht]
    \centering
    \includegraphics[width=0.45\textwidth]{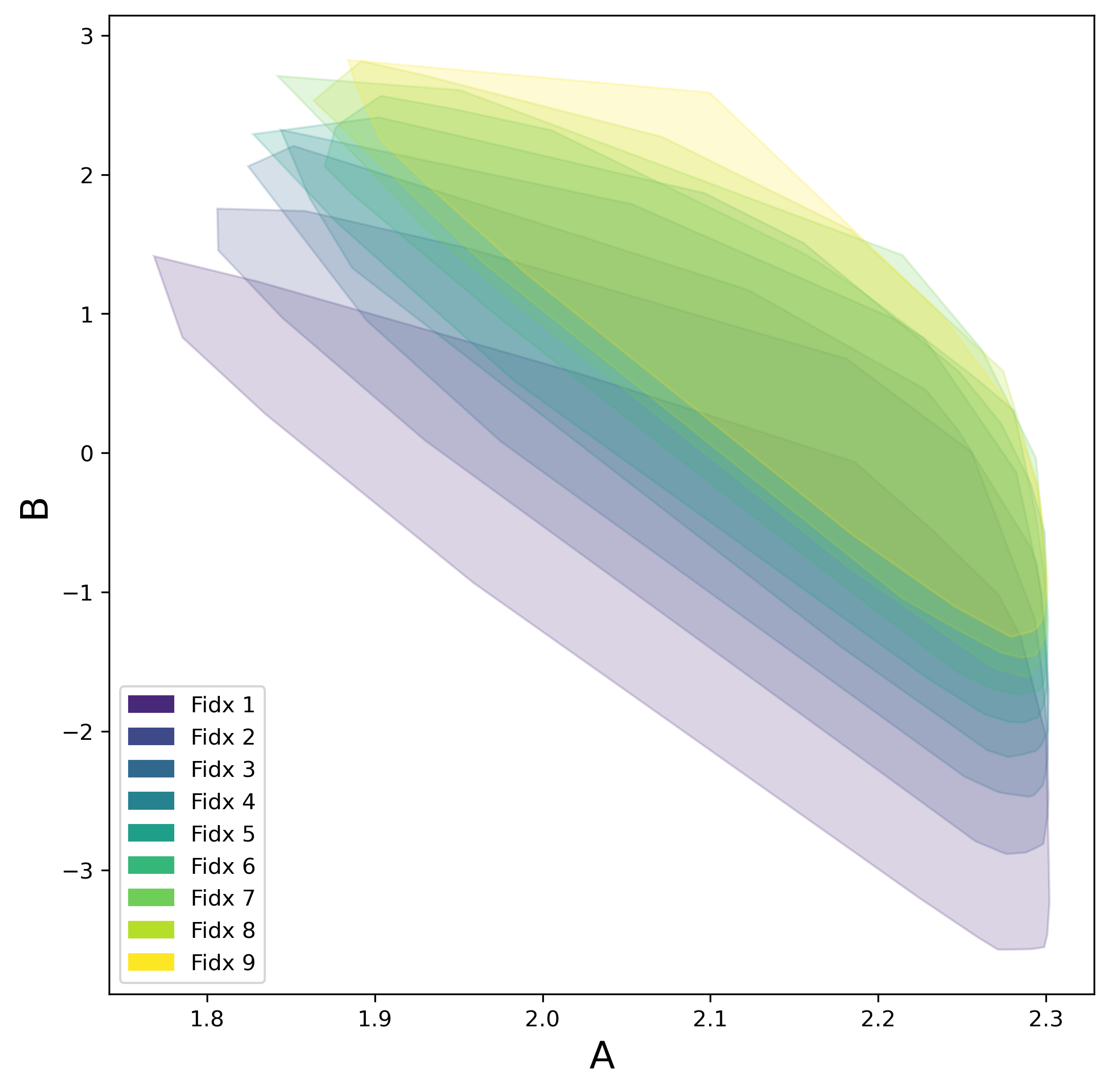}
    \caption{Convex hulls obtained from the point clouds for probably-\QFano{} toric varieties with Fano indices between~one and~nine, obtained by projecting to $\RR^2$ using the growth coefficients from~\eqref{eq:A_and_B}.}
    \label{fig:overlap}
\end{figure}

\subsection*{A cluster of high-Fano index examples}
Figure~\ref{fig:PR2_hedgehog} appears to show a cluster of high-Fano-index cases (at the top of the plot) standing apart from the remainder of the data. We now give an explanation for this high-Fano-index cluster. Figure~\ref{fig:indices_histogram} shows the frequency distribution of Fano indices in the dataset. The uptick in frequencies in the histogram in Figure~\ref{fig:indices_histogram} can be explained as follows. Consider how many ways we can write~$N$ as a sum of ten numbers between zero and seven (inclusive, and with possible repeats). This resembles a normal distribution with~$N = 35$ the most frequent case. This higher probability is due to our sampling constraints on the entries of the weight matrix: amongst those matrices that have~$a=b$ we have the highest probability of selecting one that has~$a=b=35$. Therefore, we see a misleading accumulation around those Fano indices.

In Figure~\ref{fig:overlap} we restrict the dataset to low Fano indices. For each Fano index in the range one through to nine, we plot the convex hull of the resulting point cloud. The overlap between these clusters is clear.

\section{Computational resources}\label{sec:computational_resources}
In this section we describe the computational resources required by different steps of our analysis. We will refer to a~\emph{desktop PC} and an~\emph{HPC cluster}. The desktop PC has an Intel Xeon 5222 quad-core processor, 64GB RAM, and an NVIDIA RTX A2000 12 GB GPU; note however that all CPU jobs on the desktop PC ran single-core. The HPC cluster has Intel Xeon E5-2650 processors with a total of 1400 cores.

\subsection*{Data generation}
The datasets \texttt{bound\_7\_terminal} and \texttt{bound\_7\_non\_terminal} were generated using scripts for the computational algebra system Magma~\cite{BosmaCannonPlayoust1997}, running on the HPC cluster in parallel over 1400 cores for eight days, with 2GB of memory per core. Deduplication of the dataset was performed on the desktop PC and took approximately eight hours.

\subsection*{Hyperparameter tuning}
This was carried out on the desktop PC, using the GPU. Each experiment ran on average for two minutes, for a total run time of 200 minutes for 100 experiments.

\subsection*{Model training}
This was carried out using the desktop PC, using the GPU. Training on 5M balanced samples for 150 epochs took four hours.

\subsection*{Model evaluation}
The model evaluation was carried out using the desktop PC, using the GPU. Evaluation took approximately ten minutes.

\subsection*{Further data generation}
The dataset \texttt{terminal\_dim8\_probable} was generated by running Python scripts on the HPC cluster in parallel over 120 cores for one hour, with 16GB of memory per core. Deduplication of the dataset was performed on the desktop PC and took approximately one hour.

\section{Training for weights with a larger bound}
In~\S\ref{sec:limitations} we highlighted that the trained neural network does not perform well out of sample. Therefore, it is natural to ask whether the neural network is approximating an actual general mathematical statement, or if its performance is the result of some `finite size effect' due to the choice of a particular weight bound (in our case~seven). Our intuition here is as follows. Given that the testing and training data are free of noise (they are created through exact mathematical calculation) and the neural network classifier is so accurate, we believe that the classifier is indeed approximating a precise, general mathematical statement. However, the poor out-of-sample performance makes it unclear~\emph{what kind of mathematical statement} the network is picking up. The statement could be about weight matrices with entries of arbitrary size, or could be about weight matrices with small entries (mathematically, this would be a statement about Fano varieties with terminal singularities of bounded index). In the first case the out-of-sample performance drop-off would happen because the network is approximating the true statement in a way that does not generalise to higher weight bounds; this is a common phenomenon when developing and using neural network models. In the second case the out-of-sample performance drop-off would happen because of the underlying mathematical statement that the classifier approximates. 

To probe this further, we repeated the same experiments as in the main text on a dataset of weight matrices with weights bounded by a larger constant,~ten. We generated a new dataset of size 20~million, balanced between terminal and non-terminal examples, where the entries of each weight matrix are bounded by ten. The data generation steps were the same as described in~\S\ref{sec:data_generation}, except that the terminality check was now carried out using the new algorithm discussed in~\S\ref{sec:algorithm} (and proved correct in~\S\ref{sec:proof}). We remark that the increased speed of the new algorithm allowed us to generate double the amount of data of the original dataset. 

We used a fully-connected feed-forward neural network with the same architecture as the original neural network from the paper. This architecture is recalled in Table~\ref{tab:network_architecture_supp}. Again, the network was trained on the features given by flattening the weight matrices, which where standardised by translating the mean to zero and rescaling the variance to one. It was trained using binary cross-entropy as loss function, stochastic mini-batch gradient descent optimiser and using early-stopping, for a maximum of 150~epochs and with learning rate reduction on plateaux.

\begin{figure}
    \centering
    \begin{subfigure}[b]{0.293\textwidth}
        \centering
        \includegraphics[width=\textwidth]{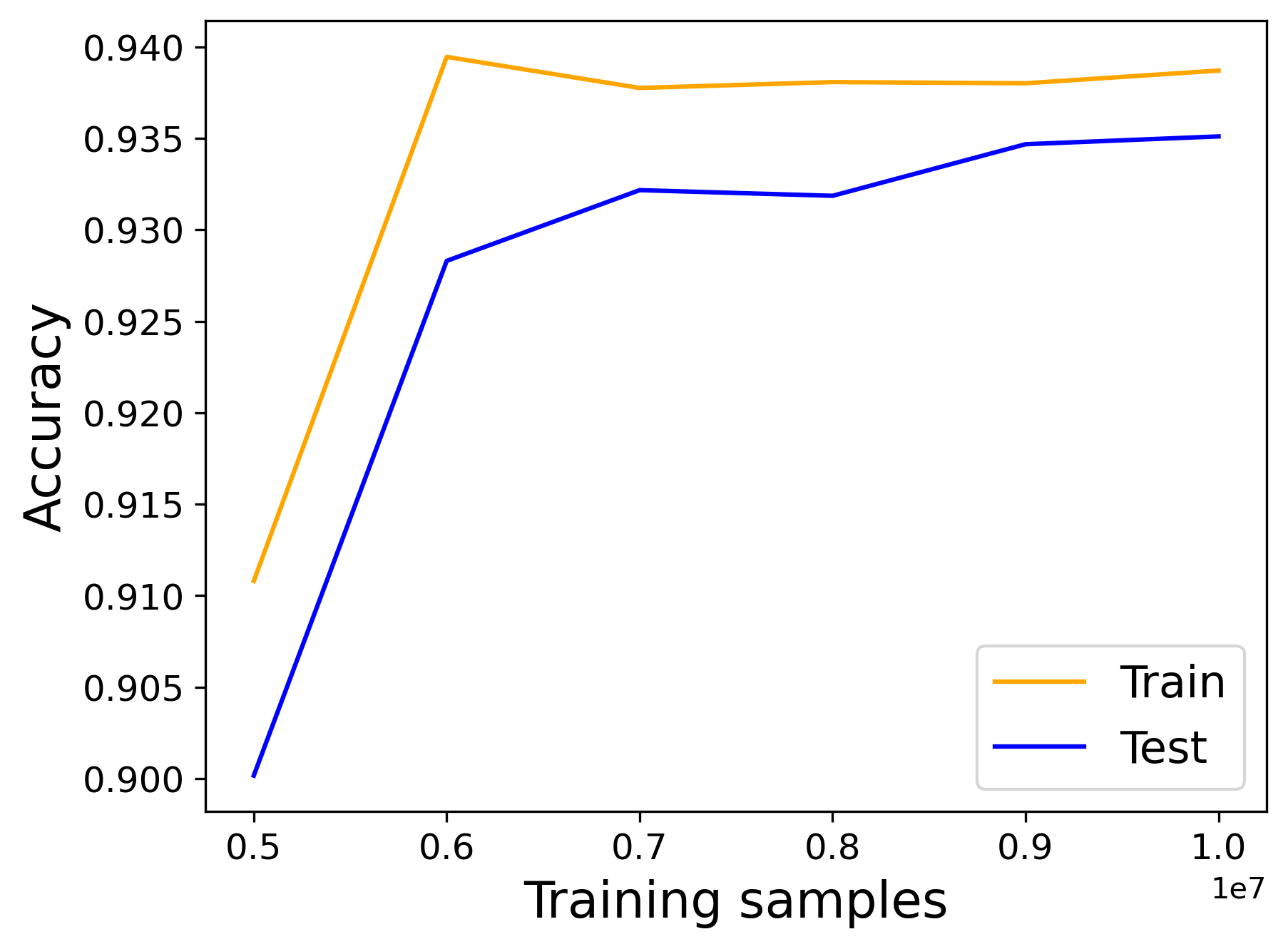}
        \caption{}\label{fig:learning_curve_train_size_bound10}
    \end{subfigure}
    \begin{subfigure}[b]{0.32\textwidth}
        \centering
        \includegraphics[width=\textwidth]{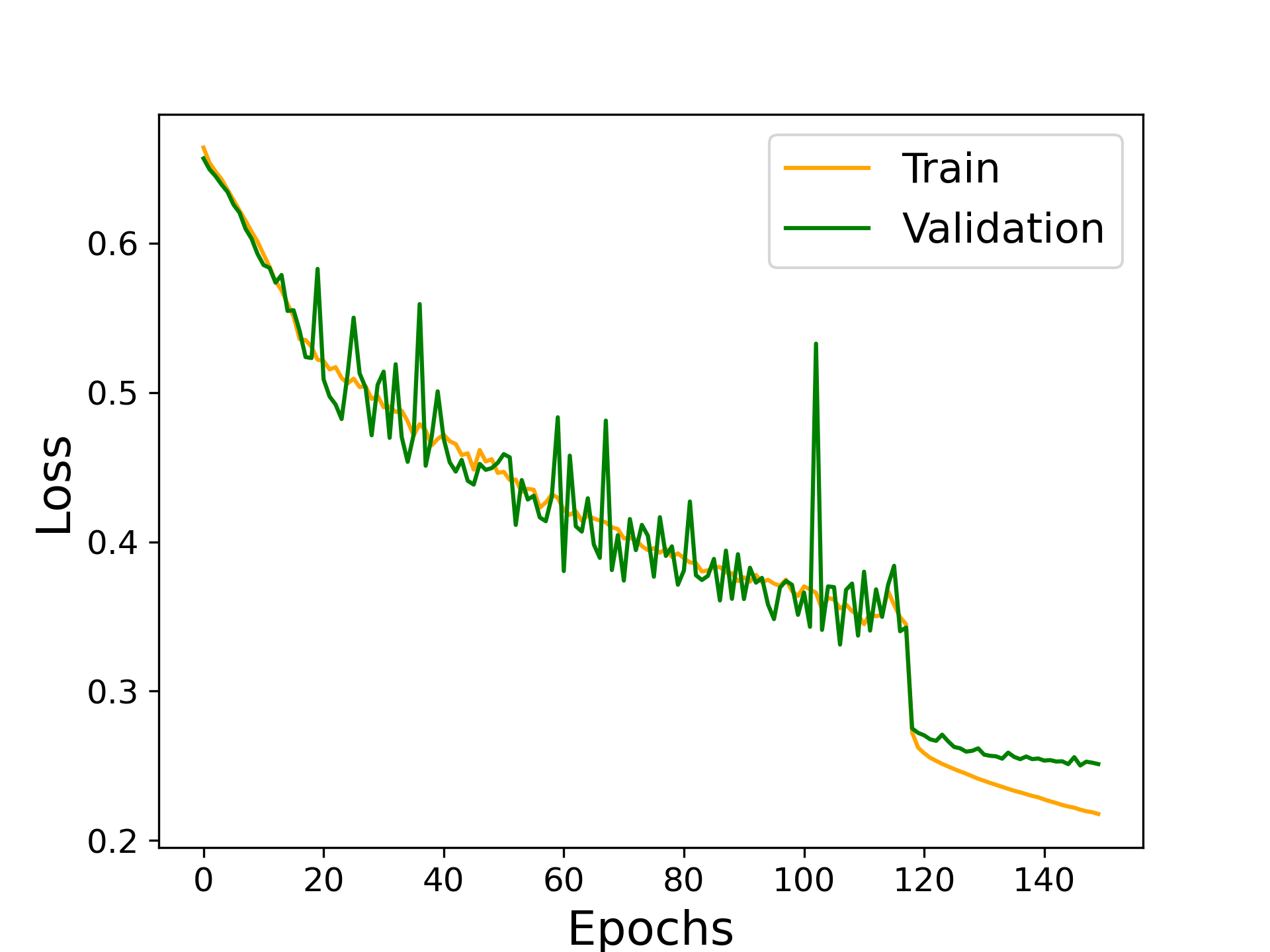}
        \caption{}\label{fig:learning_curve_terminality_small_bound10}
    \end{subfigure}
    \begin{subfigure}[b]{0.32\textwidth}
        \centering
        \includegraphics[width=\textwidth]{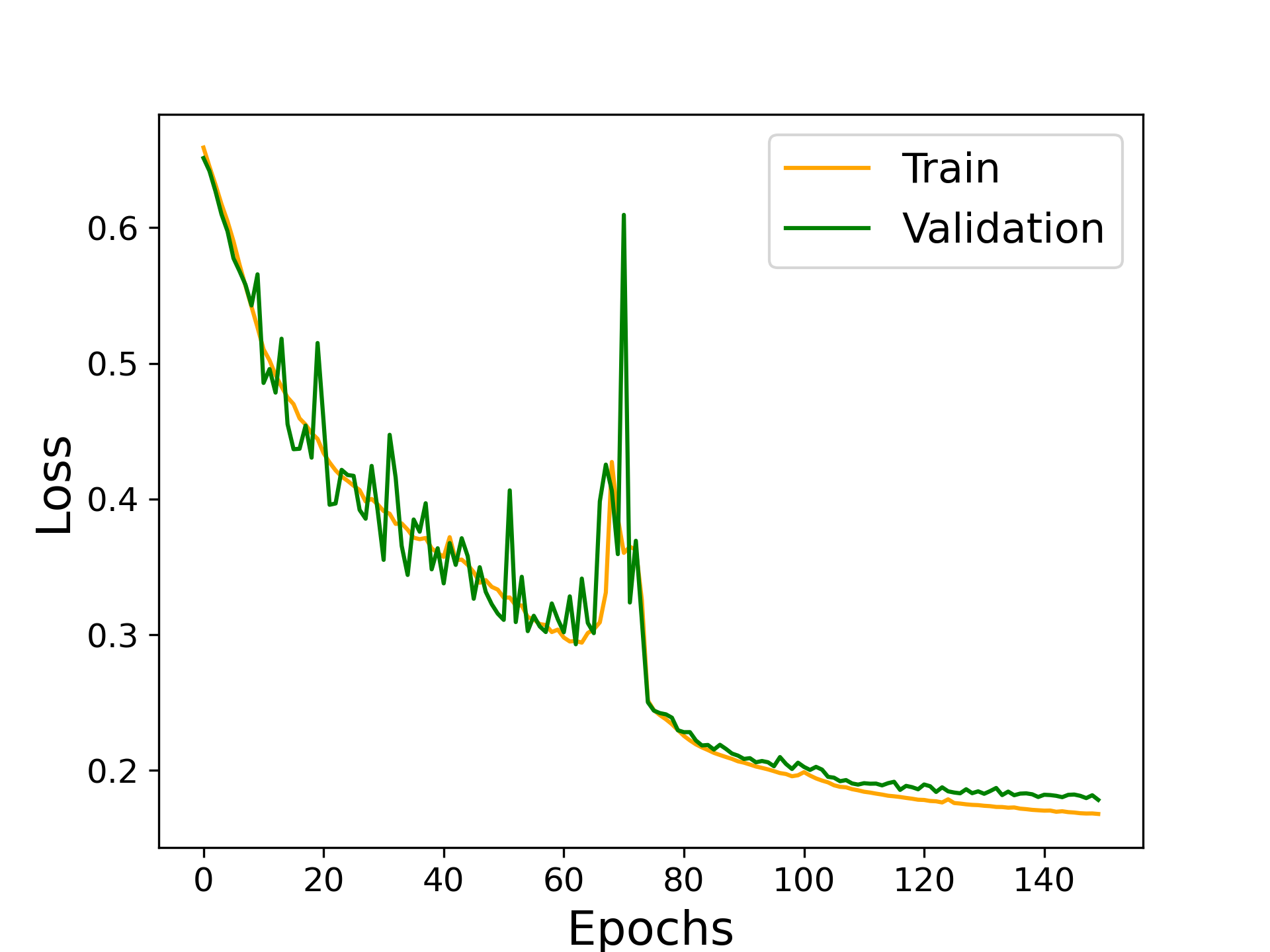}
        \caption{}\label{fig:learning_curve_terminality_big_bound10}
    \end{subfigure}
       \caption{(a) Accuracy for different train-test splits; (b) epochs against loss for the network trained on 5M samples; (c) epochs against loss for the network trained on 10M samples.}
\end{figure}

Training on 5M samples (using 80\% for training and 10\% for validation) and testing on the remaining data (15M samples) produced an accuracy of 90\% -- see Figure~\ref{fig:learning_curve_terminality_small_bound10} for the loss learning curve. This performance is worse than that achieved for the same training sample size for weight bound~seven, potentially indicating that the condition approximated by the network is harder to capture. Training on a larger sample of size 10M (again using 80\% for training and 10\% for validation) and testing on the remaining data (10M samples) improves the accuracy to 94\% -- see Figure~\ref{fig:learning_curve_terminality_big_bound10} for the loss learning curve. The training and validation accuracies for intermediate training sizes are shown in Figure~\ref{fig:learning_curve_train_size_bound10}. 

We were able to recover a high accuracy for this new dataset. However, this was only possible by using a larger training sample size, which hints at the increased difficulty of the task. Moreover, Figure~\ref{fig:learning_curve_train_size_bound10} suggests that increasing the size of the training set further is unlikely to improve the accuracy. Being able to train a high-accuracy neural network for a larger weights bound supports the hypothesis that the neural network is approximating a general mathematical statement but in a way that does not generalise well to higher bounds. However, it is too early to exclude the hypothesis that the network might be capturing a mathematical statement that needs weight matrices with small entries. Similar studies with even higher bounds would add confidence here and, if the network is indeed approximating a statement about weight matrices with small weights, experiments of this type should also be able to deduce what the cut-off bound for the weights is.

\begin{table}[t]
    \centering
    \vspace{1em}
    \begin{tabular}{rlp{0.01ex}rl}
        \toprule
        \multicolumn{1}{c}{\textbf{Hyperparameter}} &
        \multicolumn{1}{c}{\textbf{Value}} && 
        \multicolumn{1}{c}{\textbf{Hyperparameter}} & 
        \multicolumn{1}{c}{\textbf{Value}} \\
        \cmidrule{1-2} \cmidrule{4-5}
        \texttt{Layers} & $(512,768,512)$  &&
        \texttt{Momentum} & $0.99$\\
        \texttt{Batch size} & $128$ &&
        \texttt{LeakyRelu slope} & $0.01$ \\
        \texttt{Initial learning rate} & $0.01$ \\
        \bottomrule
    \end{tabular}
    \caption{Final network architecture and configuration.}\label{tab:network_architecture_supp}
\end{table}    

\begin{figure}
    \centering
    \includegraphics[width=\textwidth]{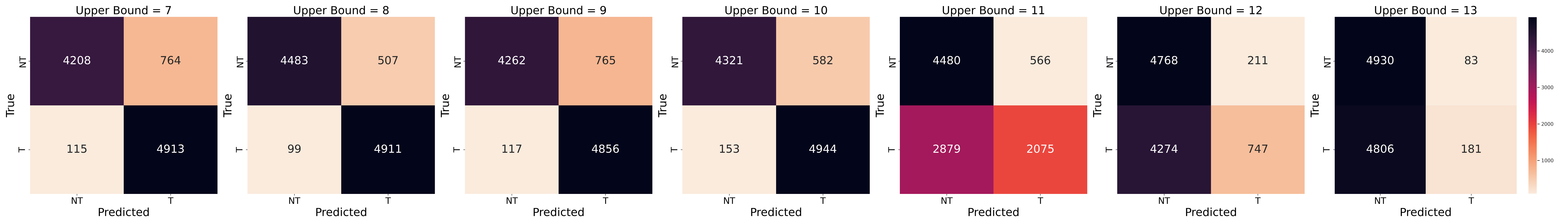}
    \caption{Confusion matrices for the neural network classifier on in-sample and out-of-sample data. In each case a balanced set of 10\,000 random examples was tested.}
       \label{fig:combined_confusion_supp}
\end{figure}

\section{Proof of Proposition~\ref{general_result}}\label{sec:proof}
In this section we prove Proposition~\ref{general_result}. This is the main ingredient in the new algorithm to check terminality. Recall from the discussion above that~$X$ determines a convex polytope~$P$ with vertices~$e_1,\ldots,e_N \in \ZZ^{N-2}$, and that 
\begin{equation*}
\begin{split}
        a_1 e_1 +\cdots +a_{N} e_{N} &= 0   \\
        b_1e_1 +\cdots +b_{N}e_{N} &= 0
\end{split}
\end{equation*}
where the~$a_i$ and~$b_j$ are entries in the weight matrix~\eqref{eq:weight_matrix_supp}. The same argument applied to the equivalent weight matrix
\[
\begin{bmatrix}
    b_i/g_i & -a_i/g_i \\
    A_i & B_i
\end{bmatrix}    
\begin{bmatrix}\label{hello}
    a_1 & \cdots & a_{N} \\
    b_1 & \cdots & b_{N}
\end{bmatrix}  
\]
gives barycentric co-ordinates for the origin and~$e_i$ in terms of the remaining vertices of~$\Delta_i$:
\begin{align*}
    \alpha^1_i e_1 + \cdots + \alpha^{i-1}_i e_{i-1} +\alpha^{i+1}_{i} e_{i+1} + \cdots + \alpha^N_{i} e_{N} & = 0 \\
    \beta^1_i e_1 + \cdots + \beta^{i-1}_{i} e_{i-1} +\beta^{i+1}_{i} e_{i+1} + \cdots + \beta^N_{i} e_{N} & = g_ie_i
\end{align*} 

Fix~$i \in \{1,2,\ldots,N\}$. Define~$u \colon \QQ^{N-1} \to \QQ$ by~$u(x_1,\ldots,x_{N-1}) = x_1 + \cdots + x_{N-1}$, and let~$\Psi$ denote the lattice
\[
    \{v \in \mathcal{Z} \mid u(v) =1 \}
\]
where~$\mathcal{Z}$ is the span over~$\ZZ$ of the standard basis~$E_1, \dots , E_{N-1}$ for~$\QQ^{N-1}$ together with
\begin{align*}
    \frac{1}{f_i} (\alpha_1^2, \dots, \hat{\alpha}_i^i, \ldots, \alpha_1^N) && \text{and} &&
    \frac{1}{g_i} (\beta_1^2, \dots , \hat{\beta}_i^i, \ldots, \beta_1^N) 
\end{align*}
Here the~$\hat{\ }$ indicates that the~$i$th entry in each vector is omitted. We define~$\phi\colon\Psi \rightarrow \ZZ^{N-2}$ to be the~$\ZZ$-linear map that sends~$E_1, \dots , E_{N-1}$ to~$e_1, \ldots, \hat{e}_i, \ldots, e_{N}$ and
\begin{align*}
    \phi\left(\frac{1}{f_i} (\alpha_1^2, \dots, \hat{\alpha}_i^i, \ldots, \alpha_1^N)\right) = 0 &&
    \phi\left(\frac{1}{g_i} (\beta_1^2, \dots , \hat{\beta}_i^i, \ldots, \beta_1^N)\right) = e_i
\end{align*}
It is easy to see that~$\phi$ is well-defined and bijective. 

Consider the higher\nobreakdash-dimensional parallelepiped~$\Gamma$ in~$\mathcal{Z}$ generated by the standard basis of~$\ZZ^{N-1}$. We note that each lattice point of~$\mathcal{Z}$ in~$\Gamma$ can represented as a linear combination
\begin{align}\label{eq:representative}
    \frac{k}{f_i} (\alpha_1^2, \dots, \hat{\alpha}_i^i, \ldots, \alpha_1^N) + \frac{l}{g_i} (\beta_1^2, \dots , \hat{\beta}_i^i, \ldots, \beta_1^N) 
\end{align}
for some~$k \in \{0, 1, \dots, f_i -1\}$ and~$l \in \{ 0,1, \dots, g_i-1\}$; this representation is unique if and only if the vertices of~$\Delta_i$ span~$\ZZ^{N-2}$. Hence,~$\Delta_i$ is almost empty if and only if whenever
\begin{align}\label{eq:terminality_relation}
    \sum_{j \ne i} \left\{ k\frac{\alpha_i^j}{f_i} + l\frac{\beta_i^j}{ g_i}  \right\} = 1
\end{align}
we have that the linear combination in~\eqref{eq:representative} represents the origin. But this is the case if and only if
\[
    \left\{ k\frac{\alpha_i^j}{f_i} + l\frac{\beta_i^j}{g_i}\right\} = \left\{ \frac{\alpha_i^j}{\alpha_i}\right\}
\]
for all~$j$, since~$(k,l) = (\frac{f_i}{\alpha_i}, 0)$ represents the origin by construction. Note that the sum~\eqref{eq:terminality_relation} could include~$j=i$, since that term is an integer and its fractional part will not contribute to the sum. \qed

\bibliographystyle{plain}
\bibliography{bibliography}
\end{document}